\documentclass[journal,twoside,web]{ieeecolor}
\pdfoutput=1
\makeatletter
\let\NAT@parse\undefined
\makeatother

\usepackage{generic}
\usepackage[sort, numbers, compress]{natbib}

\usepackage{amsthm}
\usepackage{amsmath,amssymb,amsfonts}

\usepackage{algorithmic}
\usepackage{graphicx}
\usepackage{algorithm,algorithmic}
\usepackage{hyperref}
\usepackage{textcomp}

\usepackage{subcaption}
\usepackage{mathtools}
\usepackage{hyperref}
\usepackage{enumerate}
\usepackage{comment}
\usepackage{multirow}
\def\rot{\rotatebox}

\usepackage{accents}
\newcommand\munderbar[1]{\underline{#1}}
\newtheorem{thm}{Theorem}
\newtheorem{lem}{Lemma}

\newtheorem{rem}{Remark}
\newtheorem{assumption}{Assumption}

\DeclareMathOperator*{\argmin}{arg\,min}
\allowdisplaybreaks[4]

\makeatletter
\long\def\@maketablecaption#1#2{\@tablecaptionsize
  \setbox\@tempboxa\hbox{#1. #2}
  \ifdim \wd\@tempboxa >\hsize              
    \unhbox\@tempboxa\par                   
  \else                                     
    \global \@minipagefalse
    \hbox to\hsize{\hfil\box\@tempboxa\hfil}
  \fi}
\makeatother

\def\RR{\mathbb{R}}

\def\defeq{\coloneqq}
\def\X{\mathcal{X}}
\def\U{\mathcal{U}}
\def\W{\mathcal{W}}

\usepackage{parskip}
\usepackage[algo2e,ruled,vlined,linesnumbered,norelsize]{algorithm2e}
\makeatletter
\newcommand{\removelatexerror}{\let\@latex@error\@gobble}
\makeatother

\def\BibTeX{{\rm B\kern-.05em{\sc i\kern-.025em b}\kern-.08em
    T\kern-.1667em\lower.7ex\hbox{E}\kern-.125emX}}
\begin{document}


\title{A successive convexification approach for robust receding horizon control}
\author{Yana Lishkova and Mark Cannon
\thanks{Y.L. acknowledges funding from the EPSRC Doctoral
Training Partnership EP/R513295/1, project reference
2280382 }
\thanks{Yana Lishkova and Mark Cannon are with the Department of Engineering Science, University of Oxford, OX1 3PJ, UK (e-mail: yana.lishkova@eng.ox.ac.uk,mark.cannon@eng.ox.ac.uk).}}


\maketitle

\begin{abstract}
A novel robust nonlinear model predictive control strategy is proposed for systems with nonlinear dynamics and convex state and control constraints. Using a sequential convex approximation approach and a difference of convex functions representation, the scheme constructs tubes that contain predicted model trajectories, accounting for approximation errors and disturbances, and guaranteeing constraint satisfaction. An optimal control problem is solved as a sequence of convex programs. 
We develop the scheme initially in the absence of external disturbances and show that the proposed nominal approach is non-conservative, with the solutions of successive convex programs converging to a locally optimal solution for the original optimal control problem. We extend the approach to the case of additive disturbances using a novel strategy for selecting linearization points. As a result we formulate a robust receding horizon strategy with guarantees of recursive feasibility and cloesd loop stability. 
\end{abstract}


\section{Introduction}
\label{sec:introduction}
\IEEEPARstart{M}{odel} Predictive Control (MPC), also known as receding horizon control, is an optimization-based control strategy applicable to systems with state and input constraints, model uncertainty and external disturbances. The method relies on solving an optimization problem at each discrete time step. 
The feasibility of this problem and its solution within a single time step are  
crucial for implementing the resulting control law \cite{smpc}. For nonlinear systems, the optimal control problem is in general a nonlinear program (NLP), which is not guaranteed to be solvable within a fixed time interval.

To address this problem various approaches have been developed. For example, dynamic programming can determine globally optimal control policies, however they require heavy computation, which restricts their implementation in real-time applications \cite{smpc}. Another method, known as sequential quadratic programming (SQP) \cite{10.1007/978-3-0348-8407-5_14}, focuses on iterative approximation approaches that approximate NLPs using sequences of quadratic programming (QP) problems. Despite the existence of efficient methods for solving QPs, these schemes can lead to high computational demand and present difficulties when accounting for approximation errors and determining the appropriate step-sizes for each iteration to ensure convergence.

An alternative approach replaces the system model with a sequence of models obtained from linear approximations around trajectories of the controlled system (e.g.~\cite{rosen66,mao18,OrigPaper,LKC02}). 
The seminal work~\cite{rosen66} relaxes the nonlinear model dynamics into a set of inequalities, which are subsequently linearized to obtain a sequence of convex optimization problems. Under convexity assumptions on the nonlinear dynamics and the state and input constraints, the solutions of this sequence are shown to converge to a local optimum of the original NLP provided a set of parameters (which are introduced into the objective function to enforce satisfaction of the original model dynamics) are chosen appropriately. However, the required values of cost parameters are not usually known, and the method converges only asymptotically to a feasible solution of the NLP.

More recently, \cite{mao18} proposed using Jacobian model linearization to construct a sequence of convex optimization problems. The authors use trust regions to account for the effects of linearization errors and are thus able to demonstrate convergence to a locally optimal solution of the NLP. The disadvantages of the approach are that it relies on heuristics for determining whether model approximations are sufficiently accurate, and that the iterates are not guaranteed to be feasible for the NLP, so the iteration cannot be terminated prior to convergence.

In~\cite{OrigPaper}, a tube (or sequence of sets) containing all possible trajectories of the nonlinear system is constructed by linearizing the model around trajectories and bounding linearization errors, and all trajectories within the tube are required to respect the system state and control constraints. 
Using successive linearization, the NLP is replaced by a sequence of second-order cone programs (SOCPs)~\cite{OrigPaper,LKC02}, each of which generates a feasible solution for the original NLP.
The implementation of the approach is complicated by the need to provide bounds on the linearization errors, which can be a computationally difficult task that relies on heuristics. Conservative linearization error bounds can cause slow convergence or infeasibility of the sequence of optimization problems even if the NLP itself is feasible. Furthermore, the available schemes can introduce suboptimality and do not account for external disturbances. \textcolor{black}{The tube-based approach in~\cite{foust2020optimal} can ensure convergence to a locally optimal solution of the NLP. However, it is computationally costly and does not account for external disturbances.} 

Despite the disadvantages of SQP and successive convexification approaches based on linearizing nonlinear dynamic models, the efficiency and practicality of convex NLP approximations are well-known in the context of receding horizon control. Many real-time applications of MPC therefore employ convex descriptions of the model dynamics and constraints (in aerospace applications for example, \cite{DSCF21,ZPVL18,acikmese2007convex} among many others). Moreover, successive convexification incorporating robust tube MPC methods provide new opportunities for computationally efficient receding horizon control, particularly when the model is affected by unknown disturbances.

\textcolor{black}{Approaches that consider external disturbances and model inaccuracies are reviewed in \cite{VQDCH17,KSMA21,sasfi2023robust}. The tube-based strategies of~\cite{MKWF11,YMCA13,VQDCH17,KSMA21,sasfi2023robust} address both continuous and discrete systems with various cost functions and disturbance realisations. However, these approaches use invariant sets to define tubes, rely on heuristics to determine parameters or approximate the NLP formulation conservatively. Alternatively, \cite{LSH19} proposes a two-tier approach that simultaneously optimizes tubes and trajectories online for feedback linearizable or minimum phase systems. On the other hand, \cite{MD21} optimizes a linear feedback law, alternating between a Riccati Recursion and a nominal optimal control problem, with approximate robustness due to unaccounted linearization errors in the nonlinear constraints.}

\textcolor{black}{The approach presented here exploits the convexity properties of systems with continuous dynamics
and convex state, control and disturbance constraints. We use tubes to bound the effects of model approximations and disturbances, and exploit tight bounds on linearization errors to obtain non-conservative solutions for the disturbance-free case. 
Although the min-max approach presented relies on enumerating polytopic tube vertices, which can computationally costly for high-order systems, we  discuss strategies to ensure linear growth in the number of variables and constraints with the state dimension.}


\textcolor{black}{We assume that the system discrete model 
has a convex difference (DC) form, i.e.~it can be represented as the difference of a pair of convex functions. This is not necessarily restrictive, as 
any continuous function can be approximated to any desired precision as the difference of two convex functions~\cite{horst1999dc}.} Differences of convex functions have been used to develop techniques such as DC programming \cite{horst1999dc}, the convex-concave procedure~\cite{sriperumbudur2009convergence,LB16} \textcolor{black}{and novel control invariant set computation methods for nonlinear systems \cite{fiacchini2010computation}}, which have proven numerical advantages when applied to nonconvex problems (\cite{ahmadi18,GB20}). Using DC decomposition we derive a tube-based receding horizon controller with robustness to disturbances and guarantees of recursive feasibility and stability. 

\textcolor{black}{To this end,} the first part of this paper considers the nominal
(disturbance-free) case in which the system model is known exactly. The proposed approach solves a sequence of convex programs optimizing control variables and tube cross-sections simultaneously, similar to~\cite{DSC22}. Our contributions are a novel initialisation strategy and a proof that the resulting optimization procedure is non-conservative in that its iterations converge to an (at least locally) optimal point of the NLP.

In the disturbance-free case, 
feasibility is guaranteed if at least one trajectory of the system satisfies the constraints of the problem. 
%
However, when unknown disturbances are present, model trajectories under all realisations of the disturbance must be accounted for, and this presents non-trivial additional difficulties for choosing linearisation points and guaranteeing recursive feasiblity.
%
%
The second part of this paper addresses these problems by considering additive disturbances explicitly. As in the disturbance-free case, the proposed approach approximates the NLP with a sequence of convex programs, and a tube containing all possible trajectories is optimized online. Based on the implementation of a novel linearization technique, we prove recursive feasibility and a form of quadratic stability for the resulting closed-loop system.




The contributions of this paper are as follows:
\begin{itemize}
    \item A novel tube-based receding horizon strategy with guaranteed stability and recursive feasibility for nonlinear systems with additive disturbances.
    \item Analysis of convergence of the proposed nominal tube-based approach to a local minimum of the NLP.
    \item A novel initialisation method for both the nominal case and the case in which additive disturbances are present.
\end{itemize}

\section{Successive linearization MPC}

Consider a discrete time nonlinear system with state $x_i\in \RR^{n_x}$ and control input $u_i \in \RR^{n_u}$  at  the $i$th time step, where
\begin{equation}
    x_{i+1} = f(x_i,u_i) , \quad  i =0,1,\ldots
    \label{eqn:nonlinearmodel}
\end{equation} 
and $f$ is differentiable for all $(x_i,u_i)$ in the operating region $\mathcal{X}\times\mathcal{U}$.  We  pose  the  problem  of  optimally controlling the system trajectories $\mathbf{x}:\{x_0,x_1,...\}$, $\mathbf{u}:\{u_0,u_1,...\}$ to track a given state and input reference trajectory $\mathbf{x^r}:\{x^r_0,x^r_1,...\}$, $\mathbf{u^r}:\{u^r_0,u^r_1,...\}$ where $x^r_{i+1} = f(x^r_i, u^r_i)$ for all $i \geq 0$ subject to a quadratic objective 
\begin{align}
\refstepcounter{equation}
 &\sum_{i=0}^{\infty} \Big( \left\| x_{ i}-x^r_{i} \right\|^2_Q+\left\| u_{ i}-u^r_{i}\right\|^2_R \Big) \label{eqn:ncost}
\end{align}  
  and initial constraints $x_0 = x_{\textrm{init}}$ and  inequality constraints 
\begin{align} \label{eqn:stateandinputconstr}
\refstepcounter{equation}
 &   x_i \in \mathcal{X} \subset {\rm I\!R}^{n_x}, \;\;  u_i \in \mathcal{U} \subset {\rm I\!R}^{n_u}.
\end{align}  
Here $\mathcal{X}$, $\mathcal{U}$ are compact polytopic sets, $x_{\textrm{init}}$ is a known initial condition and $\| x \|^2_Q \defeq x^TQx$. 
To solve this problem we define the following receding horizon optimal control problem (RHOCP):
\begin{subequations} \label{eqn:nMPC}
\begin{alignat}{2}
&\textrm{nMPC:}\quad
& (\mathbf{u}^*_i,\mathbf{x}^*_i) & \defeq 
\argmin_{\mathbf{u}_i,\mathbf{x}_i}
J_{\textrm{nMPC}}(\mathbf{u}_i, \mathbf{x}_i) 
\label{eqn:NMPC1}\\
&\text{\makebox[0pt][l]{subject to, for $k = 0, \ldots , N-1:$}} 
\nonumber\\
&&& x_{k \vert i} \in \mathcal{X}  
\label{eqn:constr2} \\
&&&  u_{k \vert i} \in \mathcal{U}  
\label{eqn:constr3} \\
&&&  x_{k+1 \vert i} = f(x_{k \vert i},u_{k \vert i})  
\label{eqn:constr4}\\
&&&  x_{0 \vert i} = x_{i}  
\label{eqn:constr1}\\ 
&&& 
x_{N \vert i} \in \mathcal{X}_N \subseteq  \mathcal{X} ,
\label{eqn:constr5}
\end{alignat}
with
\begin{equation}\label{eqn:nmpccost}%
\begin{aligned}
& J_{\textrm{nMPC}}(\mathbf{u}_i,\mathbf{x}_i) \defeq
\| x_{N\vert i} - x^r_{N+i} \|^2_{Q_N} 
\\
&\qquad + \sum_{k=0}^{N-1}  \bigl( \| x_{k \vert i} - x^r_{k+i} \|^2_Q 
+  \| u_{k \vert i} - u^r_{k+i} \|^2_R\bigr) ,
\end{aligned}
\end{equation}
\end{subequations}
where $\smash{\mathbf{u}_i\defeq\{u_{0 \vert i },\ldots,u_{N-1 \vert i }\}}$ and
$\smash{\mathbf{x}_i\defeq\{x_{0 \vert i },\ldots,x_{N \vert i }\}}$. Here $x_{k \vert i}$, $u_{k \vert i}$ denote the state and control input at time $k+i$ predicted at time $i$. The matrices $Q,R,Q_N$ and the terminal set $\mathcal{X}_N$ are subject to the following assumptions: 
\begin{assumption} \label{assumpt:matrices}
(a)~$x_i$ is known at each time $i$.\\
(b)~$x^r_i\in \mathcal{X}_N\subseteq\X$ and $u_i^r\in\mathcal{U}$, $\forall i\geq N$.\\
(c)~$f(x,\kappa_i(x))\in\X_N$, $\forall x\in\X_N$, $\forall i\geq N$, where $\kappa_i(x) \defeq K(x - x^r_i) + u^r_i$ satisfies $\kappa_i(x) \in \mathcal{U}$, $\forall x\in\mathcal{X}_N$, $\forall i\geq N$.\\
(d)~$Q, R, Q_N$ are positive definite matrices, and \eqref{eqn:nonlinearmodel} with 
$u_{i} = K(x_{i} - x^r_{i}) + u^r_{i}$ 
satisfies, $\forall x_{i}\in \mathcal{X}_N$, $\forall i\geq N$,
\[
\|x_i - x_i^r\|_{Q_N}^2 - \|x_{i+1} - x_{i+1}^r \|_{Q_N}^2 
\geq 
\| x_i - x^r_i \|_Q^2 + \| u_i - u^r_i \|_R^2 .
\]
\end{assumption}
In receding horizon control, at each discrete time step $i$ the nMPC RHOCP is solved to obtain an optimal control sequence $\mathbf{u}_i^{*}:\{u^*_{0 \vert i },u^*_{1 \vert i },\ldots,u^*_{N-1 \vert i }\}$ and the first element of this sequence is used to define a receding control law $u_i=u^*_{0 \vert i }$. At the next time step a new RHOCP is formulated using the measurement (or estimate) of $x_{i+1}$ 
and the process is repeated. To guarantee stability of the control strategy it is essential that a feasible solution can be obtained within each discrete time step. However, in its form in \eqref{eqn:NMPC1}-\eqref{eqn:nmpccost}, the nMPC RHOCP is a nonlinear program and thus in general no guarantees can be given that either the optimal solution or a suboptimal feasible point is determinable within one discrete time step.
To address these problems Section~\ref{section:nominalsection} proposes a successive linearization approach to reformulate the problem as a sequence of convex programs. Compared to other such approaches, the proposed method bounds linearization errors as functions of the predicted state tube that is optimized online. This removes the need to compute linearization error bounds prior to solving the RHOCP and ensures that a fixed point of the successive approximation scheme corresponds to an optimal point of the nMPC RHOCP. As shown in Section~\ref{section:proofs}, it also provides guarantees of convergence to a fixed point, recursive feasibility and stability of the resulting MPC law.   

\section{Proposed control law} \label{section:nominalsection}
To simplify explanation we first strengthen the assumptions on the system model and later discuss how these assumptions can be relaxed to allow the proposed approach to be applied to systems with dynamics that are neither convex or concave. 

\begin{assumption} \label{assumpt:two}
The function $f: \mathcal{X}\times\mathcal{U} \to \RR^{n_x}$ in \eqref{eqn:nonlinearmodel} is differentiable and componentwise convex.
\end{assumption}%

 Given state and control sequences $\mathbf{x}^0_i:\{x_{0 \vert i}^0,x_{ 1 \vert i }^0,\ldots,x_{ N \vert i }^0\}$ and $\mathbf{u}^0_i:\{u_{ 0 \vert i }^0,u_{ 1 \vert i }^0,\ldots,u_{ N-1 \vert i }^{0}\}$, which we refer to as a seed trajectory, 
 %
 the model \eqref{eqn:nonlinearmodel} can be equivalently expressed as
\begin{gather}
    x_{k+1 \vert i} = x^0_{k+1 \vert i}+ A_{k \vert i}  s_{k \vert i}+  B_{k \vert i}v_{k \vert i} + e(x_{k \vert i},u_{k \vert i}) \label{eqn:Taylor}
\\
s_{k \vert i} = x_{k\vert i}-x^0_{k\vert i} , \;\;\;\; v_{k\vert i} =u_{k\vert i}-u^0_{k\vert i}\;\;  \label{eqn:pertrmodel}
\\
   A_{k \vert i} = \frac{\partial f}{\partial x_{k \vert i}} \vert_{(x_{k \vert i}^0,u_{k \vert i}^0)}, 
    B_{k \vert i} = \frac{\partial f}{\partial u_{k \vert i}} \vert_{(x_{k \vert i}^0,u_{k \vert i}^0)} \; 
    \label{eqn:statematrices}
\end{gather}
and $e_{k|i}$ is the Jacobian linearization error
\begin{align} \label{eqn:error}
\!\!\! e_{k|i}  = f(x_{k|i},\!u_{k|i})- f(x_{k|i}^0,\!u_{k|i}^0) - A_{k|i}s_{k|i}-B_{k|i}v_{k|i} .
\end{align}
 \begin{assumption} \label{assumpt:seed}
    The initial seed trajectory $(\mathbf{x}^0,\mathbf{u}^0)$ is feasible (but potentially suboptimal) for constraints  \eqref{eqn:constr1}-\eqref{eqn:constr5}.
\end{assumption}

Anticipating the explicit treatment of disturbances in Section~\ref{section:robustsection}, we assume a feedback control law:
\begin{equation}\label{eqn:ueqn}
u_{k \vert i} = u^0_{k \vert i}+v_{k \vert i}=u_{k \vert i}^0+Ks_{k \vert i} + c_{k \vert i}, 
\end{equation}
for a given feedback gain $K$,
and consider the sequence $\mathbf{c}_i\defeq \{c_{0 \vert i}, \ldots , c_{N-1 \vert i} \}$ as an optimization variable instead of $\mathbf{u}_i$.

Assumption~\ref{assumpt:two} requires the components of $f$
to be convex and hence the components of the linearization error $e_{k|i}$ in~\eqref{eqn:error}, 
denoted $[e_{k \vert i}(x_{k \vert i},u_{k \vert i})]_j$, $j= 1,\ldots,n_x$, are convex functions of $(x_{k\vert i}, u_{k \vert i})\in(\mathcal{X},\, \mathcal{U})$. 
Since $[e_{k|i}(x_{k \vert i},u_{k \vert i})]_j$ is convex and the linearization point $(x^0_{k \vert i},u^0_{k \vert i})$ lies in $(\X,\,\U)$ we have
\begin{equation}\label{eqn:min_lin_error_bound}
\min_{(x_{k \vert i}, u_{k|i}) \in (\X,\, \U)} [e_{k \vert i}(x_{k \vert i},u_{k \vert i})]_j = 0 ,
\end{equation}
and moreover a constraint of the form
\[
[e_{k \vert i}(s_{k \vert i}+x^0_{k \vert i}, Ks_{k \vert i} + c_{k \vert i}+ u^0_{k \vert i})]_j \leq r 
\]
is convex in $s_{k|i}$, $c_{k|i}$ and $r$.

Let $\mathbf{S}_i=\{S_{0 \vert i},\ldots, S_{N \vert i}\}$ denote a tube bounding the state perturbations $s_{k|i}$ in~\eqref{eqn:Taylor}-\eqref{eqn:statematrices}, where each set $S_{k \vert i}$ is defined in terms of optimization variables $\smash{\munderbar{s}_{k \vert i}}$, $\bar{s}_{k \vert i} \in \mathbb{R}^{n_x}$ as
\begin{align} \label{eqn:tubes}
    S_{k \vert i} = \{ s : \munderbar{s}_{k \vert i} \leq s \leq \bar{s}_{k \vert i} \} .
\end{align}
Then at time $i$ a convex RHOCP can be formulated as
\begin{subequations} \label{eqn:nominalSOCP}
\begin{alignat}{2}
&\textrm{cMPC:}\hspace{1mm}
& & (\mathbf{c}_i^*,\mathbf{S}_i^*) \defeq \argmin_{\mathbf{c}_i,\,\mathbf{S}_i}
\, J_{\textrm{cMPC}}(\mathbf{c}_i,\mathbf{S}_i,\mathbf{x}^0_i, \mathbf{u}^0_i)
\hspace{50mm}
\nonumber\\
&\text{\makebox[0pt][l]{subject to, for $k = 0, \ldots , N-1:$}} 
\nonumber\\
&&& S_{k \vert i} \oplus \{x^0_{k \vert i} \}  \subseteq  \mathcal{X} \label{eqn:cMPCRHOCP2}\\
&&& KS_{k|i} \oplus \{c_{k \vert i} + u^0_{k \vert i}\} \subseteq  \mathcal{U} 
\label{eqn:cMPCRHOCP3}\\
&&&  [\bar{s}_{k+1 \vert i}]_j \geq [A_{k \vert i}+B_{k \vert i}K]_j s + [B_{k \vert i}]_j c_{k \vert i}  \nonumber \\
&&& \qquad  \qquad +[e_{k \vert i}(s+x^0_{k \vert i},Ks+c_{k \vert i}+u^0_{k \vert i})]_j \nonumber \\
&&& [\munderbar{s}_{k+1 \vert i}]_j \leq [A_{k \vert i}+B_{k \vert i}K]_j s + [B_{k \vert i}]_j c_{k \vert i}  \nonumber  \\
&&&  \qquad  \qquad \ \ \text{for all } s\in S_{k|i}, \, j = 1,\ldots,n_x
\label{eqn:cMPCRHOCP4}  \\
&&& S_{0 \vert i} = \{0\}
\label{eqn:cMPCRHOCP1}\\
&&& S_{N \vert i} \oplus \{ x_{ N \vert i}^0\} \subseteq \mathcal{X}_N 
\label{eqn:cMPCRHOCP5}
\end{alignat}
with
\begin{equation}\label{eqn:cMPCRHOCP0}
\begin{aligned}
&  J_{\textrm{cMPC}}(\mathbf{c}_i,\mathbf{S}_i,\mathbf{x}^0_i, \mathbf{u}^0_i) \defeq \!\max_{s_{N \vert i} \in S_{N \vert i}}\! \| x^0_{N\vert i} + s_{N\vert i}- x^r_{N+i}\|^2_{Q_N} 
\\
&\quad +\sum_{k=0}^{N-1}
\max_{s_{k \vert i} \in S_{k \vert i}} \Bigl[
\| x_{k \vert i}^0+s_{k \vert i}-x^r_{k+i}\|^2_Q
\\
&\quad\qquad\qquad\qquad+
\| u_{k \vert i}^0+Ks_{k \vert i}+c_{k|i}- u^r_{k+i}\|^2_R\Bigr] ,
\end{aligned}
\end{equation}
\end{subequations}
where $[z]_j$, $[A]_j$ denote the $j$th element and $j$th row of vector $z$ and matrix $A$. 
The inequalities in \eqref{eqn:cMPCRHOCP4} can be imposed for all $s\in \smash{S_{k|i}}$ by replacing $s$ with the $\smash{2^{n_x}}$ vertices of $\smash{S_{k|i}}$, which depend linearly on  $\smash{\munderbar{s}_{k \vert i}}$, $\bar{s}_{k \vert i}$. The minimization of the cost~\eqref{eqn:cMPCRHOCP0} can be performed via an equivalent epigraph form with a linear objective subject to convex constraints on ${\bf c}_i$ and $\{\smash{\munderbar{s}_{k \vert i}}$, $\bar{s}_{k \vert i} , \ k= 0,\ldots,N\}$
the vertices of $\smash{S_{k|i}}$. The cMPC RHOCP is convex and can be solved using, for example, interior-point methods or large-scale solvers such as the alternating direction method of multipliers (ADMM)~\cite{Reading72,boyd11}.

\textcolor{black}{The feedback gain $K$ can be chosen equal to the terminal feedback gain $K$, which is stabilizing within $\mathcal{X}_N$ due to Assumption~\ref{assumpt:matrices}(d). The computation of $K$, $\mathcal{X}_N$, and $Q_N$ is in general non-trivial for nonlinear systems. A scheme to compute them while maximizing the volume of the terminal set using semidefinite programming can be found in \cite{lishkova2022multirate,smpc}.  
}

The cMPC strategy is summarized in Algorithm~\ref{alg:cap}. At each discrete time step $i$, the algorithm first solves the cMPC RHOCP to obtain $\mathbf{c}_{i}^*$ given a feasible seed $(\mathbf{x}^0_i,\mathbf{u}^0_i)$ (lines 3-4), then computes the predicted state and control sequences $(\mathbf{x}_i ^\ast,\mathbf{u}_i ^\ast)$ using the model dynamics~\eqref{eqn:nonlinearmodel} with $\mathbf{c}_i=\mathbf{c}_i^\ast$ and $x_{0|i}^\ast=x_i$ (lines 5-6), then updates the seed trajectory $(\mathbf{x}_i^0, \mathbf{u}_i^0)$ (line 7) and repeats this iteration until either $\|\mathbf{c}_i^\ast\| < \mathit{tol}$, where $\mathit{tol}$ is a specified threshold, or the limit $\mathit{maxiters}$ on the number of iterations is reached.

{\setlength{\algomargin}{1.4em}
\removelatexerror
\begin{algorithm2e}
\SetKwInOut{Input}{Offline}
\SetKwInOut{Output}{Online}
\Input{Determine $K$, $Q_N$, $\mathcal{X}_N$, and an initial seed trajectory $(\mathbf{x}_0^0, \mathbf{u}_0^0)$ satisfying \eqref{eqn:constr2}-\eqref{eqn:constr5} with $x_{0} = x_{\mathrm{init}}$}
\Output{At each discrete time step $i=0,1,\ldots$:}
Obtain $x_i$ and set the iteration counter $n \gets 1$\;
\While{$n \leq \mathit{maxiters}$ and $\|\mathbf{c}_i^* \| \geq \mathit{tol}$}{
Compute $A_{k \vert i}$, $B_{k \vert i}$ for $k=0,\ldots,N-1$ using \eqref{eqn:statematrices} and $(\mathbf{x}_i^0, \mathbf{u}_i^0)$\;
Solve cMPC RHOCP to obtain $(\mathbf{c}_i^\ast,\mathbf{S}_i^\ast)$\;
\For{$k=0,\ldots,N-1$}{
Compute $u_{k \vert i}^\ast= u^0_{k \vert i} +c^*_{k \vert i}+K(x_{k \vert i}^\ast-x^0_{k \vert i})$  and $x_{k+1\vert i}^\ast = f(x_{k \vert i}^\ast,u_{k \vert i}^\ast)$ with $x_{0 \vert i}^\ast = x_{i}$\;
}
Update the seed trajectory: $(\mathbf{x}_i^0, \mathbf{u}_i^0) \gets (\mathbf{x}_i ^\ast, \mathbf{u}_i ^\ast )$\;
 $n \gets n + 1$\;
}
Implement $u_i \gets u_{0 \vert i}^\ast$\;
Set $u_{N|i}^\ast = K(x_{N\vert i}^\ast - x^r_{N+i}) + u^r_{N+i}$ and
$\mathbf{u}_{i+1}^0 \gets \{ u_{1 \vert i}^\ast,\ldots,u_{N-1 \vert i}^\ast, u_{N \vert i}^\ast\}$  
$\mathbf{x}_{i+1}^0 \gets \{ x_{1 \vert i}^\ast,\ldots,x_{N \vert i}^\ast, f(x_{N|i}^\ast,u_{N|i}^\ast) \}$\;
\caption{cMPC}\label{alg:cap}
\end{algorithm2e}}

\begin{rem}
A time-varying feedback gain (such as the linear-quadratic optimal controller for the linear time-varying model computed on line~3 of Algorithm~1) may be used instead of a fixed feedback gain $K$ in \eqref{eqn:ueqn}. This is likely to give better numerical conditioning and disturbance rejection since the fixed gain $K$ is only required by Assumption~\ref{assumpt:matrices} to be locally stabilizing within the terminal set $\X_N$. 
\end{rem}

Compared to other successive linearization approaches~\cite{mao18,OrigPaper,LKC02}, the cMPC strategy does not require prior knowledge or online calculation of sets bounding linearization errors, and it avoids conservative approximate linear constraints by using tight bounds on the approximation error in \eqref{eqn:cMPCRHOCP4}. In particular, the first inequality in \eqref{eqn:cMPCRHOCP4} is equivalent to $\bar{s}_{k+1 \vert i,j} \geq [f(s_{k \vert i}+x^0_{k \vert i}, Ks_{k \vert i}+c_{k \vert i}+u^0_{k \vert i})]_j-[f(x^0_{k \vert i},u^0_{k \vert i})]_j$, and the second is a tight bound due to \eqref{eqn:min_lin_error_bound}.
As shown in Section \ref{section:proofs}, this ensures that the tube computed by the cMPC iteration converges to a single trajectory, and hence a fixed point of the cMPC iteration is a local minimum for the nMPC RHOCP. In addition, the maximum number of iterations can be chosen as any value satisfying $\mathit{maxiters}\geq 1$ without affecting the recursive feasibility and stability guarantees discussed in Section~\ref{section:proofs}.

A procedure for obtaining an initial seed trajectory is suggested in \cite{DSC22}. This assumes knowledge of an initial trajectory satisfying constraints \eqref{eqn:constr2}-\eqref{eqn:constr4} and the initial condition \eqref{eqn:constr1}, but not necessarily the terminal constraint \eqref{eqn:constr5}. The method updates this seed via a sequence of convex optimizations to obtain a feasible trajectory for the cMPC RHOCP.

We propose an alternative method of computing offline an initial seed trajectory $(\smash{\mathbf{x}_0^0},\smash{\mathbf{u}_0^0})$ for Algorithm~\ref{alg:cap} that does not require knowledge of a feasible trajectory with initial condition $x_{\mathrm{init}}$. We define an offline optimization of $(\mathbf{c}_0,\mathbf{S}_0)$ subject to the constraints \eqref{eqn:cMPCRHOCP2}-\eqref{eqn:cMPCRHOCP4}, \eqref{eqn:cMPCRHOCP5}, and $ S_{0|0}=\{s\}$, where $s$ is an additional optimization variable, and we minimize the distance of $\smash{x_{0|0}^0}+s$ from $x_{\mathrm{init}}$ by setting the objective to $\|x_{\mathrm{init}} - \smash{x_{0|0}^0}-s\|$. 
Initially we set $\mathbf{x}^0_0\defeq \{x^r_N , \ldots,x^r_{2N}\} $ and $\mathbf{u}^0_0\defeq\{u^r_N , \ldots , u^r_{2N-1} \} $, since by Assumption~\ref{assumpt:matrices} this choice ensures that constraints \eqref{eqn:cMPCRHOCP2}-\eqref{eqn:cMPCRHOCP4}, \eqref{eqn:cMPCRHOCP5} and $S_{0|0}=\{s\}$ are feasible with $s=0$. To obtain an initial seed for Algorithm~\ref{alg:cap}, we iteratively solve this optimization and update $(\smash{\mathbf{x}_0^0},\smash{\mathbf{u}_0^0})$ as in lines 2-8 of Algorithm~\ref{alg:cap}, with $x_{0|0}^\ast= x_0$ in line 6 replaced by $x_{0|0}^\ast= x_{0|0}^0+s^\ast$, until either the optimal objective is equal to zero, or the optimal objective converges to a non-zero value, in which case the constraints should (if possible) be relaxed, for example by increasing $N$ or enlarging $\mathcal{X}_N$.

\begin{rem}\label{rem:nominalDC}
Assumption \ref{assumpt:two} can be relaxed to allow functions $f$ which are DC, i.e. they can be written as the difference of convex functions: $f(x,u) = g(x,u) - h(x,u)$ where $g,h: \mathcal{X}\times \mathcal{U} \to \RR^{n_x}$  are differentiable and componentwise convex and the sets $\mathcal{X},\mathcal{U}$ are convex. In this case the error function has the following properties 
\begin{align*} 
&e_{k|i} = e^g_{k|i}-e^{h}_{k|i} \\
&e^g_{k|i}  = g(x_{k|i},u_{k|i})- g(x_{k|i}^0,u_{k|i}^0) - A^g_{k|i}s_{k|i}-B^g_{k|i}v_{k|i} \\
&e^h_{k|i}  = h(x_{k|i},u_{k|i})- h(x_{k|i}^0,u_{k|i}^0) - A^h_{k|i}s_{k|i}-B^h_{k|i}v_{k|i} \\
&\min_{(x_{k \vert i}, u_{k|i}) \in (\X,\, \U)} [e^g_{k \vert i}(x_{k \vert i},u_{k \vert i})]_j=0  \\&\max_{(x_{k \vert i}, u_{k|i}) \in (\X,\, \U)} -[e^{h}_{k \vert i}(x_{k \vert i},u_{k \vert i})]_j = 0 \notag
\end{align*}
where $e^g_{k \vert i},e^h_{k \vert i}$ are componentwise convex and 
\begin{align*}
   &A^g_{k \vert i} = \frac{\partial g}{\partial x_{k \vert i}} \vert_{(x_{k \vert i}^0,u_{k \vert i}^0)}, \;
    B^g_{k \vert i} = \frac{\partial g}{\partial u_{k \vert i}} \vert_{(x_{k \vert i}^0,u_{k \vert i}^0)}  \\
       &A^h_{k \vert i} = \frac{\partial h}{\partial x_{k \vert i}} \vert_{(x_{k \vert i}^0,u_{k \vert i}^0)}, 
    \;B^h_{k \vert i} = \frac{\partial h}{\partial u_{k \vert i}} \vert_{(x_{k \vert i}^0,u_{k \vert i}^0)} 
\end{align*}
Thus Algorithm \ref{alg:cap} can be used with constraints \eqref{eqn:cMPCRHOCP4} in step 4 replaced, for all  $s\in S_{k|i}$, $j = 1,\ldots,n_x$, by
\begin{align*}
  [\bar{s}_{k+1 \vert i}]_j &\geq [A_{k \vert i}+B_{k \vert i}K]_j s + [B_{k \vert i}]_j c_{k \vert i} \\
&\quad +[e^g_{k \vert i}(s+x^0_{k \vert i},Ks+c_{k \vert i}+u^0_{k \vert i})]_j \\
[\munderbar{s}_{k+1 \vert i}]_j &\leq [A_{k \vert i}+B_{k \vert i}K]_j s + [B_{k \vert i}]_j c_{k \vert i} \\ 
&\quad -[e^h_{k \vert i}(s+x^0_{k \vert i},Ks+c_{k \vert i}+u^0_{k \vert i})]_j .
\end{align*}
\end{rem}

\begin{rem} \label{rem:sizenominal} If the cMPC RHOCP \eqref{eqn:nominalSOCP} is implemented by replacing $s$ with the $\smash{2^{n_x}}$ vertices of $\smash{S_{k|i}}$ (which depend linearly on  $\smash{\munderbar{s}_{k \vert i}}$, $\bar{s}_{k \vert i}$), then the number of constraints in \eqref{eqn:cMPCRHOCP4} will grow exponentially with the state dimension $n_x$.
However, this exponential growth can be avoided by using homothetic tubes (see e.g.~\cite{smpc}) to bound the state perturbations $s_{k|i}$. In this formulation the tube description in \eqref{eqn:tubes} is replaced with 
    \[
        S_{k \vert i} = \{z^0_{k\vert i}\} \oplus \alpha_{k|i} \mathrm{co}\{z^1,...,z^{n_v}\} ,
    \]
    where $\mathrm{co}$ denotes the convex hull. Here  $z^0_{k\vert i}\in\mathbb{R}^{n_x}$ and $\alpha_{k|i} \in\mathbb{R}$ are optimization variables in the RHOCP and the vertices $\{z^1,...,z^{n_v}\}$ are fixed by the designer. This reduces the number of optimization variables at each time step from $2n_x$ to $n_x+1$ and reduces the number of constraints in \eqref{eqn:cMPCRHOCP4} to $2m n_x$. Here $n_v$ can be as small as $n_x+1$ if $S_{k|i}$ is chosen to be a simplex,
    allowing the number of constraints in the cMPC RHOCP to depend quadratically on the state dimension $n_x$. \textcolor{black}{This is discussed in more detail in Remark~\ref{rem:simplex_tube} in Section \ref{section:robustsection}.}%
\end{rem}

\section{Feasibility and stability guarantees} \label{section:proofs}
In this section we demonstrate that Algorithm~\ref{alg:cap} is recursively feasible and that it provides a feasible solution for the nMPC RHOCP. 
We further show that the optimal cMPC RHOCP objective value is an upper bound on the cost of the nMPC RHOCP, and that the optimal cost is non-increasing at successive iterations (lines 2-8 of Algorithm~\ref{alg:cap}) and at successive time steps. 
Using these properties, we show that the iteration converges (if the iteration limit $\mathit{maxiters}$ is sufficiently large) to a solution satisfying the first-order optimality conditions of the nMPC RHOCP, and that (for any value of $\mathit{maxiters}$) the control law of Algorithm~\ref{alg:cap} is asymptotically stable. 

\begin{thm}[Recursive feasibility] \label{thm:recrsfeas}
If Assumptions \ref{assumpt:matrices} and \ref{assumpt:two} hold and the initial seed trajectory $(\mathbf{x}^0_0,\mathbf{u}^0_0)$ satisfies \eqref{eqn:constr2}-\eqref{eqn:constr5}, then the cMPC RHOCP is feasible at each iteration of Algorithm~\ref{alg:cap}, at all times $i \geq 0$.
\end{thm}
    
\textbf{Proof.}\hspace{1ex}
Let $\mathbf{c}_i = 0$ and $\mathbf{S}_{i} = \{\{0\},\ldots,\{0\}\} $, then the constraints of the cMPC RHOCP in line 4 of Alg.~\ref{alg:cap} reduce to
$\{\smash{x_{k | i}^0}\} \subseteq  \mathcal{X}$, $\smash{u_{k| i}^0} \in \mathcal{U}$ and $\{ \smash{x_{ N| i}^0} \} \subseteq \mathcal{X}_N$.
These conditions are satisfied if $(\mathbf{x}_i,\mathbf{u}_i) =(\mathbf{x}_i^0,\mathbf{u}_i^0)$ satisfies \eqref{eqn:constr2}-\eqref{eqn:constr5}, which is true by assumption for the initial seed trajectory $(\mathbf{x}_0^0,\mathbf{u}_0^0)$, and is true for the seed trajectories  $(\mathbf{x}_i^0,\mathbf{u}_i^0)$ and  $(\mathbf{x}_{i+1}^0,\mathbf{u}_{i+1}^0)$ in lines 7 and 10 of Alg.~\ref{alg:cap} due to constraints \eqref{eqn:cMPCRHOCP2}-\eqref{eqn:cMPCRHOCP5} and the properties of $\mathcal{X}_N$ (Assumption~\ref{assumpt:matrices}(b) and (c)). \qed 

We next consider how the optimal cMPC RHOCP cost relates to the cost associated with the seed trajectory at successive iterations of Algorithm~\ref{alg:cap}, and hence prove that the iteration converges to a solution that is at least locally optimal for the nMPC RHOCP. 
At the $i$th time step and $n$th iteration of Algorithm~\ref{alg:cap}, let $(\smash{\mathbf{x}^{0,n}_i},\smash{\mathbf{u}^{0,n}_i})$ denote the seed trajectory in line 3, let $(\mathbf{c}^{*,n}_i,\mathbf{S}^{*,n}_i)$ denote the optimal solution of the cMPC RHOCP in line 4 and define
$\smash{J_{\textrm{cMPC},i}^{*,n}} \defeq J_{\textrm{cMPC}} (\mathbf{c}^{*,n}_i,\mathbf{S}^{*,n}_i, \smash{\mathbf{x}^{0,n}_i}, \smash{\mathbf{u}^{0,n}_i} ) $.
  
\begin{lem} \label{lem:cost_sequence}
For all $i\geq 0$ and for all $n\geq 1$, we have 
$
J_{nMPC}( \mathbf{x}_i^{0,n+1}, \mathbf{u}_i^{0,n+1} ) \leq J_{\textrm{cMPC},i}^{*,n} \leq
J_{nMPC}( \mathbf{x}_i^{0,n}, \mathbf{u}_i^{0,n} )
$.
\end{lem}

\textbf{Proof.}\hspace{1ex}
Let $\mathbf{c}_i =0$ and $\mathbf{S}_i=\{\{0\},\ldots,\{0\}\}$, then
\[
 J_{\textrm{cMPC}} (0, \{\{0\},\ldots,\{0\}\}, \mathbf{x}^{0,n}_{i}\! , \mathbf{u}^{0,n}_{i} ) 
=
 J_{\textrm{nMPC}} (\mathbf{x}^{0,n}_i\!,\mathbf{u}^{0,n}_i)  ,
\]
and since this is a feasible (possibly suboptimal) choice, the optimal cMPC cost in line 4 of Algorithm \ref{alg:cap} satisfies
\[
 J_{\textrm{cMPC}} (\mathbf{c}_i^\ast, \mathbf{S}_i^\ast, \mathbf{x}^{0,n}_{i}, \mathbf{u}^{0,n}_{i} ) 
\leq
J_{\textrm{nMPC}} (\mathbf{x}^{0,n}_i,\mathbf{u}^{0,n}_i) .
\]
Furthermore, the updated seed in lines 5-7 satisfies
\[
J_{\textrm{nMPC}} (\mathbf{x}^{0,n+1}_i,\mathbf{u}^{0,n+1}_i) 
\leq
 J_{\textrm{cMPC}} (\mathbf{c}_i^\ast, \mathbf{S}_i^\ast, \mathbf{x}^{0,n}_{i}, \mathbf{u}^{0,n}_{i} ) 
\]
since the tube $\mathbf{S}_i^\ast$ necessarily contains $\mathbf{x}_i^{0,n+1}$.
\qed 

\begin{thm} \label{thm:convergence}
For all $i\geq 0$, the iteration on lines 2-8 of Algorithm~\ref{alg:cap} converges to a seed trajectory $(\smash{\mathbf{x}^{0,n}_i},\smash{\mathbf{u}^{0,n}_i})$ such that $\mathbf{c}_i^{\ast,n}=0$ and $\mathbf{S}_i^{\ast,n} = \{\{0\}, \ldots, \{0\}\}$ is an optimal solution of the cMPC RHOCP in the limit as $n\to\infty$.
\end{thm}

\textbf{Proof.}\hspace{1ex}
From Lemma \ref{lem:cost_sequence} we have $J_{\textrm{cMPC},i}^{*,n+1} \leq J_{\textrm{cMPC},i}^{*,n}$ for all $n$, so $\smash[t]{J_{\textrm{cMPC},i}^{*,n}}\geq 0$ implies $\smash[t]{J_{\textrm{cMPC},i}^{*,n}}-\smash[t]{J_{\textrm{cMPC},i}^{*,n+1}}\to 0$ and hence ${J_{\textrm{cMPC},i}^{*,n}} - {J_{\textrm{nMPC}}(\mathbf{x}_i^{0,n},\mathbf{u}_i^{0,n})}  {\to 0}$ as $n\to\infty$. Therefore $(\mathbf{c}_i^{\ast,n},\mathbf{S}_i^{\ast,n}) = (0,\{\{0\}, \ldots, \{0\}\})$ is an optimal solution of the cMPC RHOCP in the limit as $n\to\infty$.
\qed
    
\begin{rem}
\label{rem:convergence}
The optimal solution of the cMPC RHOCP for $\mathbf{c}^\ast_i$ is unique since $R$ is positive definite. Therefore Theorem~\ref{thm:convergence} implies $\mathbf{c}^{\ast,n}_i\to 0$ as $n\to \infty$, justifying the termination criterion in line 2 of Alg.~\ref{alg:cap}. Moreover, even if the optimal $\mathbf{S}^{\ast}_i$ is non-unique, we can choose $\mathbf{S}^{\ast,n}_i=0$ if $\mathbf{c}^{\ast,n}_i =0$ and thus ensure that the iteration converges to a tube containing a single trajectory.
\end{rem} 
    
\begin{thm} \label{thm:KKT}
The iteration of Algorithm 1 converges to a local minimum of the nMPC RHOCP \eqref{eqn:constr2}-\eqref{eqn:nmpccost}.
\end{thm} 

A proof of Theorem~\ref{thm:KKT} can be found in the Appendix. 
    
\begin{rem} \label{rem:optima}
The converse of Theorem~\ref{thm:KKT} is also true: a local minimum of the nMPC problem is a local minimum of the cMPC RHOCP \eqref{eqn:cMPCRHOCP2}-\eqref{eqn:cMPCRHOCP0} due to Theorem~\ref{thm:convergence}. 
\end{rem}

\begin{thm} \label{thm:assymstability}
Given Assumptions~\ref{assumpt:matrices} and \ref{assumpt:two}, the control law of Algorithm 1 ensures that $x-x^r = 0$ is an asymptotically stable equilibrium of the system~\eqref{eqn:nonlinearmodel} \footnote{More precisely, $x^\delta = x-x^r = 0$ is an asymptotically stable equilibrium of the system $x^\delta_{i+1} = f(x^\delta_i + x_i^r , u_i) - f(x_i^r,u_i^r)$.}, with region of attraction consisting of the initial conditions $x_0$ of \eqref{eqn:nonlinearmodel} for which there exists a control sequence such that  $u_{k}\in \mathcal{U}$ and $x_{k}\in \mathcal{X}$ for $k=0,\ldots,N-1$, and $x_{N} \in \mathcal{X}_N$.
\end{thm} 
    
\textbf{Proof.}\hspace{1ex}%
Assumption~\ref{assumpt:matrices}(d) and the update policy in line~10 of Algorithm~\ref{alg:cap} imply that the cost associated with the seed trajectory at the first iteration at time $i+1$ satisfies
\begin{align*}
J_{\textrm{nMPC}} (\mathbf{x}_{i+1}^{0,1}, \mathbf{u}_{i+1}^{0,1}) 
&\leq 
 J_{\textrm{nMPC}}  (\mathbf{x}_{i}^{0,n_i}, \mathbf{u}_{i}^{0,n_i}) 
\\
&\quad - (\| x_{i} - x^r_{i} \|_Q^2 + \| u_{i} - u_{i}^r \|_R^2 ) ,
\end{align*}
where $n_i$ is the final iteration count at time $i$.
Hence 
\begin{align*}
J_{\textrm{nMPC}}(\mathbf{x}_{i+1}^{0,n_{i+1}}\!\!,\mathbf{u}_{i+1}^{0,n_{i+1}}\!)
&\leq 
J_{\textrm{nMPC}}(\mathbf{x}_{i}^{0,n_i},\mathbf{u}_{i}^{0,n_i}) 
\\
&\quad - (\| x_{i} - x^r_{i} \|_Q^2 + \| u_{i} - u_{i}^r \|_R^2 ) 
\end{align*} 
due to Lemma~\ref{lem:cost_sequence},
and, since $Q,R,Q_N$ are positive definite, $J_{\textrm{nMPC}}(\mathbf{x}_{i+1}^{0,n_i},\mathbf{u}_{i+1}^{0,n_i})$ is therefore a Lyapunov function demonstrating stability of $x-x^r=0$. Summing both sides of this inequality over $i=0,\ldots,t$ we obtain
\[
\sum_{i=0}^\infty (\| x_{i} - x_{i}^r \|_Q^2 + \|u_{i}  - u_{i} ^r\|_R^2) \leq
J_{\textrm{nMPC}}({\bf x}_0^{0,n_1}, {\bf u}_0^{0,n_1} )
\]
and it follows that $(x_{ i} , u_{i} ) \to ( x_{ i}^r, u_i^r)$ as $i\to \infty$. \qed 

Theorem \ref{thm:KKT} and Remark \ref{rem:optima} show that the cMPC RHOCP formulation allows Algorithm~\ref{alg:cap} to effectively approximate the optimal control law by solving a sequence of convex programs. 

\section{Robust MPC with additive disturbance} \label{section:robustsection}
This section considers uncertain systems of the form
\begin{align}\label{eqn:robustsystemmodel}
x_{i+1} = f(x_i,u_i) + w_i, \qquad w_i\in\W
\end{align}
where $w_i$ is an unknown disturbance input and $\W\subset\mathbb{R}^{n_x}$ is a known compact set. We consider an optimal tracking control problem similar to \eqref{eqn:nMPC}, with the model~\eqref {eqn:robustsystemmodel} replacing \eqref{eqn:nonlinearmodel}. 
In this setting we retain Assumption \ref{assumpt:two} 
(we discuss later how this can be relaxed to allow nonconvex dynamics), 
and we replace Assumption \ref{assumpt:matrices} with the following assumption.

\begin{assumption}\label{assumpt:robust}
(a)~$x_i$ is known at each time $i$.\\
\textcolor{black}{(b)~$\mathcal{X}_N\subseteq \mathbb{R}^{2n_x}$ is known such that, 
for all $(\munderbar{x}_N,\bar{x}_N)\in \mathcal{X}_N$ and all $x\in X(\munderbar{x}_N,\bar{x}_N) := \{x: \munderbar{x}_N \leq x \leq \bar{x}_N\}$
we have: 
(i)~$x^r_i\in X(\munderbar{x}_N,\bar{x}_N)$ and $u_i^r\in\mathcal{U}$ $\forall i\geq N$,
(ii) $X(\munderbar{x}_N,\bar{x}_N) \subseteq \X$ and $(\munderbar{x}_N^+,\bar{x}_N^+)\in\mathcal{X}_N$ where, $\forall i\geq N$, 
\begin{align*}
[\bar{x}_{N}^+]_j &\geq [f(x,\kappa_i(x))]_{j} +[\bar{w}]_j
\\
[\munderbar{x}_{N}^+]_j &\leq [f(x,\kappa_i(x))]_{j} + [\munderbar{w}]_j  , \quad j = 1,\ldots,n_x ,
\end{align*}
and
(iii) $\kappa_i(x) \defeq K(x - x^r_i) + u^r_i \in \mathcal{U}$, $\forall i\geq N$.}\\
(c)~$Q, R, Q_N \succ 0$ and \eqref{eqn:robustsystemmodel} with  $u_{i} = K(x_{i} - x^r_{i}) + u^r_{i}$ satisfies, for all $i\geq N$, for some $\beta \geq 0$ and some $(\munderbar{x}_N,\bar{x}_N)\in \mathcal{X}_N$,
\begin{align*}
&\|x_i - x_i^r\|_{Q_N}^2 - \|x_{i+1} - x_{i+1}^r \|_{Q_N}^2 
\\ 
&\quad \geq 
\| x_i - x^r_i \|_Q^2 + \| u_i - u^r_i \|_R^2 - \beta, \ \ \forall x_{i}\in X(\munderbar{x}_N,\bar{x}_N) .
\end{align*}
\end{assumption}


The uncertainty in the future disturbance sequence implies that future state trajectories generated by \eqref{eqn:robustsystemmodel} for a given control sequence or control law are uncertain. Therefore, although a convex successive optimization approach can again be used to define an MPC law, the optimal tube bounding predicted model states cannot converge to a single trajectory as in the nominal case considered in Section~\ref{section:proofs}. In addition, the seed trajectory used to define a linear model approximation in the disturbance-free case must be replaced by a tube.

\subsection{Proposed robust control law}

Recursive feasibility requires that at each iteration a feasible seed tube can be constructed, guaranteeing that the state and control constraints are satisfied under all possible realizations of the disturbance input in \eqref{eqn:robustsystemmodel}. 
Given an initial state $x_{0|i}$ and perturbation sequence $\mathbf{c}_i^0$, we propose the method summarized in Algorithm~\ref{alg:cap4} for constructing a linear model approximation and a seed tube $\mathbf{X}_i^0\defeq\{\smash{X^0_{0|i}},\ldots,\smash{X^0_{N|i}}\}$, where 
\begin{equation}\label{eq:Xsetdef}
X^0_{k|i} = \{x : \ \smash{\munderbar{x}_{k|i}}\leq x \leq \smash{\bar{x}_{k|i}}\}.
\end{equation}
The linearization points, denoted $(\smash{x_{k \vert i,j}^{0},u^{0}_{k \vert i,j}})$, are computed for each component, $[f]_j$, $j=1,\ldots,n_x$, of $f$, via
\begin{equation}\label{eqn:findlinearizpoints}
     x^{0}_{k \vert i,j}= \argmin_{x \in X^0_{k\vert i} } \; [f(x,Kx +c^{0}_{k \vert i})]_j 
\end{equation}
and $\smash{u^{0}_{k \vert i,j}}  = K\smash{x^{0}_{k \vert i,j}} + \smash{c^{0}_{k \vert i}}$. We show in Section~\ref{section:RobustFeasStabGuar} that this choice ensures recursive feasibility. The seed tube is computed using elementwise bounds, for $j=1,\ldots,n_x$,
\begin{alignat}{1}
[\bar{x}_{k+1 \vert i} ]_j &=  \max_{x\in X^0_{k|i}} 
[f(x,Kx +c^{0}_{k \vert i})]_j 
+[\bar{w}]_j
\label{eqn:findxoverb}
\\
[\munderbar{x}_{k+1 \vert i} ]_j &= [f(x^{0}_{k \vert i,j},u^{0}_{k \vert i,j})]_j 
+[\munderbar{w}]_j
\label{eqn:findxunderb}
\end{alignat}
where
\[
[\bar{w}]_j = \smash{\max_{w\in\W}[w]_j}, \qquad 
[\munderbar{w}]_j = \smash{\min_{w\in\W}[w]_j} .
\]
Note that the minimization in \eqref{eqn:findlinearizpoints} is a convex problem, and the maximization in \eqref{eqn:findxoverb} can be performed as a maximization over the vertices of $\smash{X_{k|i}^0}$, each of which is a known linear function of $\smash{\munderbar{x}_{k|i}},\smash{\bar{x}_{k|i}}$.

The linear model approximation computed using Algorithm~\ref{alg:cap4} can be used to formulate the following robust RHOCP.
\begin{subequations}
\begin{alignat}{2}
&\text{\makebox[0pt][l]{$\textrm{RcMPC}:
\hspace{1mm} (\mathbf{c}_i^*,\mathbf{X}_i^*)  \defeq 
\argmin_{\substack{\;\;\;\;\;\;\;\;\mathbf{c}_i, \mathbf{X}_i}} J_{\textrm{RcMPC}}(\mathbf{c}_i,\mathbf{X}_i,\mathbf{c}^0_i)$}}
& &\hspace{100mm}
\nonumber\\
&\text{\makebox[0pt][l]{subject to, for $k = 0, \ldots , N-1:$}} 
\nonumber\\
&&& X_{k \vert i}  \subseteq  \mathcal{X} 
\label{eqn:RcMPCRHOCP3}\\
&&& KX_{k \vert i} \oplus \{c_{k \vert i} +c^0_{k \vert i}\}  \subseteq \mathcal{U} , 
\label{eqn:RcMPCRHOCP2}\\
&&& [\bar{x}_{k+1 \vert i}]_j \geq [f(x,Kx +c_{k \vert i}+c^0_{k \vert i})]_{j} 
+[\bar{w}]_j
\nonumber\\
&&& [\munderbar{x}_{k+1 \vert i}]_j \leq [f(x^{0}_{k \vert i,j},u^{0}_{k \vert i,j})]_{j} + [A_{k \vert i} \!+\! B_{k \vert i}K]_j (x-x^{0}_{k \vert i,j}) 
\nonumber\\ 
&&& 
\qquad\qquad\quad + [B_{k \vert i}]_j c_{k \vert i}  
+ [\munderbar{w}]_j  
\nonumber\\
&&& \qquad\qquad\quad \text{for all } x \in X_{k \vert i},\ j = 1,\ldots,n_x 
\label{eqn:RcMPCRHOCP4}\\
&&& X_{0 \vert i} = \{x_i\} 
\label{eqn:RcMPCRHOCP1}\\
&&& (\munderbar{x}_{N|i},\bar{x}_{N|i}) \in \mathcal{X}_N
\label{eqn:RcMPCRHOCP5}
\end{alignat}
with 
\begin{equation}\label{eqn:RcMPCRHOCP0}%
\begin{aligned}
&
J_{\mathrm{RcMPC}}(\mathbf{c}_i,\mathbf{X}_i,\mathbf{c}_i^0) \defeq \!\max_{x_{N \vert i} \in X_{N \vert i}}\! \| x_{N\vert i} - x^r_{N+i}\|^2_{Q_N} 
\\
&\quad
+ \sum_{k=0}^{N-1} \max_{x_{k \vert i} \in X_{k \vert i}} 
\Bigl[ \| x_{k \vert i} - x_{k+i}^r \|^2_Q 
\\
&\quad\qquad\qquad\qquad
+  \| Kx_{k \vert i} +c_{k \vert i}+c^0_{k \vert i}-u^r_{k+i} \|^2_R \Bigr] .
\end{aligned}
\end{equation}
\end{subequations}
The constraints in \eqref{eqn:RcMPCRHOCP4} can be imposed for all $x\in X_{k|i}$ by replacing $x$ with the $2^{n_x}$ vertices of $X_{k|i}$, which are linear functions of the optimization variables  $\munderbar{x}_{k|i}$,  $\bar{x}_{k|i}$. Hence the RcMPC RHOCP can be expressed as a convex program similar to the cMPC RHOCP for the nominal case.

{\setlength{\algomargin}{1.2em}
\removelatexerror
\begin{algorithm2e}
\SetKwInOut{Input}{Require}
\Input{Initial state $x_{i}$ and seed $\mathbf{c}_i^0$}
Set $\bar{x}_{0 \vert i} = \munderbar{x}_{0 \vert i}= x_{i}$ so that
$X^0_{0 \vert i } \gets \{x_{i} \}$\;
\For{$k=0,\ldots,N-1$}{
\For{$j=1,\ldots,n_x$}{
Solve problem \eqref{eqn:findlinearizpoints}  for $x^{0}_{k \vert i,j}$  and $u^{0}_{k \vert i,j}  = Kx^{0}_{k \vert i,j} + c^{0}_{k \vert i}$\;
Compute 
$\displaystyle [A_{k \vert i}]_j =  \frac{\partial [f]_j}{\partial x} \bigr\rvert_{(x_{k \vert i,j}^{0},u^{0}_{k \vert i,j})}$ and $\displaystyle [B_{k \vert i}]_j =  \frac{\partial [f]_j}{\partial u} \bigr\rvert_{(x_{k \vert i,j}^{0},u^{0}_{k \vert i,j})} $\;
Solve \eqref{eqn:findxoverb}-\eqref{eqn:findxunderb} for  $[\bar{x}_{k+1 \vert i}]_j$, $[\munderbar{x}_{k+1 \vert i} ]_j$\;
}
Set $X^0_{k+1 \vert i}\gets \{x : \munderbar{x}_{k+1 \vert i}\leq x \leq \bar{x}_{k+1 \vert i} \}$\;
}
\caption{Seed tube and linearization}\label{alg:cap4}
\end{algorithm2e}}

\begin{rem}
The feedback law $Kx_{k|i}+c_{k|i}+\smash{c_{k|i}^0}$ can be interpreted as a controller of the form $K s_{k|i}+c_{k|i}+\smash{v_{k|i}^0}$, where $s_{k|i} =x_{k|i}-\smash{\hat{x}_{k|i}}$, $\smash{v^0_{k|i}} = K \smash{\hat{x}_{k|i}} + \smash{c_{k|i}^0}$, and $\smash{\hat{x}_{k|i}}$ is a point in $\smash{X_{k|i}^0}$ such as the Chebyshev centre of $\smash{X_{k|i}^0}$, for example. However, this equivalent formulation would require the additional unnecessary computation of $\smash{\hat{x}_{k|i}}$ for $k=0,\ldots,N$. 
\end{rem}

The robust cMPC strategy is presented in Algorithm~\ref{alg:cap2}. The initial seed $\smash{\mathbf{c}_0^0}$ can be computed offline similarly to the initialization of Algorithm~\ref{alg:cap}. Starting with  $x_0 \defeq \smash{x^r_N}$ and $\smash{\mathbf{c}_0^0} \defeq \{\smash{u^r_N}-K\smash{x^r_N},\ldots,\smash{u^r_{2N-1}}-K\smash{x^r_{2N-1}}\}$, the tube $\smash{\mathbf{X}_0^0}$ generated by Algorithm~\ref{alg:cap4} satisfies $\smash{X^0_{k|i}}\subseteq\X_N$ for all $k$ (by Assumption~\ref{assumpt:robust}), and hence \eqref{eqn:RcMPCRHOCP3}-\eqref{eqn:RcMPCRHOCP4}, \eqref{eqn:RcMPCRHOCP5} are feasible with $(\mathbf{c}_0,\mathbf{X}_0) = (0,\mathbf{X}_0^0)$. To generate an initial seed for Algorithm~\ref{alg:cap2}, we perform the iteration in lines 2-5 of Algorithm~\ref{alg:cap2}, with line~4 replaced by the minimization of $\|x_{\mathrm{init}}-x_0\|$ over variables $(\mathbf{c}_0,\mathbf{X}_0)$ and $x_0$ subject to \eqref{eqn:RcMPCRHOCP3}-\eqref{eqn:RcMPCRHOCP4}, $\smash{X_{0|0}^0}=\{x_0\}$ and \eqref{eqn:RcMPCRHOCP5}, with $x_0\gets x_0^\ast$, ${\mathbf{c}_0^0} \gets {\mathbf{c}_0^\ast} + {\mathbf{c}_0^0}$ in line 5. This converges either to a  feasible ${\mathbf{c}_0^0}$ with $x^\ast_0 = x_{\mathrm{init}}$, or to a non-zero value of $\|x_{\mathrm{init}}-x_0^\ast\|$ (in which case the constraints should be relaxed, e.g.\ by increasing $N$ or enlarging $\X_N$).

{\setlength{\algomargin}{1.2em}
\removelatexerror
\begin{algorithm2e}
\SetKwInOut{Input}{Offline}
\SetKwInOut{Output}{Online}
\Input{Determine $K$, $Q_N$, $\mathcal{X}_N$, and initial seed $\mathbf{c}_0^0$ with which Algorithm \ref{alg:cap4} generates a tube $\mathbf{X}_0^0$ such that $(\mathbf{c}_0,\mathbf{X}_0) = (0, \mathbf{X}_0^0)$ satisfies \ref{eqn:RcMPCRHOCP3}-\eqref{eqn:RcMPCRHOCP5} with $x_{0} = x_{\mathrm{init}}$}
\Output{At each discrete time step $i=0,1,\ldots$:}
Obtain $x_i$ and set the iteration counter  $n\gets 1$\;
\While{$n \leq \mathit{maxiters}$ and $\|\mathbf{c}_i^* \| \geq \mathit{tol}$}{
Compute $(x_{k \vert i,j}^{0},u_{k \vert i,j}^{0})$, $j=1,\ldots,n_x$ and $A_{k \vert i}$, $B_{k \vert i}$ for $k=0,\ldots,N-1$ using Alg.~\ref{alg:cap4}\;
Solve RcMPC RHOCP to obtain $(\mathbf{c}_i^*, \mathbf{X}_i^*)$\;
Update the seed: $\mathbf{c}^0_{i} \gets \mathbf{c}^*_{i} + \mathbf{c}^0_{i}$, 
$n \gets n + 1$\;
}
Implement $u_i \gets K x_i+c^0_{0 \vert i}$ and set
$\mathbf{c}_{i+1}^0 \gets \{ c^0_{1 \vert i}, \ldots, c^0_{N-1 \vert i}, u^r_{N+i} - Kx^r_{N+i}\} $
\caption{RcMPC}\label{alg:cap2}
\end{algorithm2e}}

\begin{rem}
If $\W=\{0\}$,  then $(x_{k|i,j}^{0},u_{k|i,j}^{0})=(x_{k|i,1}^{0},u_{k|i,1}^{0})$ for all $j=2,\ldots,n_x$,
so at each time step $k$ all linearization points are identical. In this case the tube $\mathbf{X}_i^\ast$ generated by Algorithm~\ref{alg:cap2} contains a single trajectory, which is identical to $\mathbf{x}_i^\ast$ at each iteration of Algorithm~\ref{alg:cap}. Thus Algorithm~\ref{alg:cap} is a special case of Algorithm~\ref{alg:cap2} when disturbances are absent.
\end{rem}

Various methods are available for computing the terminal set $\X_N$, feedback gain $K$ and matrix $Q_N$ satisfying Assumption~\ref{assumpt:robust}. 
\textcolor{black}{The set $\X_N$ can be chosen for example as the the maximum robust positively invariant (MRPI) set or as a sub-optimal (non-maximal) terminal set which satisfies Assumption~\ref{assumpt:robust}. Methods for computing  $K$, $Q_N$, and the MRPI offline, including a constraint-checking approach, can be found in~\cite{buerger2024,smpc,LORENZEN2019461,blanchini2008set}} 

\begin{rem}\label{rem:DC_dist}
Analogously to the disturbance-free case in Remark~\ref{rem:nominalDC}, Assumption \ref{assumpt:two} can be relaxed to consider DC functions $f(x,u) = g(x,u)-h(x,u)$, where $g,h: \mathcal{X}\times\mathcal{U} \to \mathbb{R}^{n_x}$ are differentiable and componentwise convex on the convex sets $\mathcal{X},\mathcal{U}$. In this case Algorithms \ref{alg:cap4} and \ref{alg:cap2} can still be used, with \eqref{eqn:findlinearizpoints} in step 4 of Algorithm \ref{alg:cap4} replaced, for $\phi = g,h$ by
\[
x^{0,\phi}_{k \vert i,j} = \argmin_{x \in X^0_{k\vert i} } \; [\phi(x,Kx +c^{0}_{k \vert i})]_j, 
\]
and with \eqref{eqn:findxoverb}, \eqref{eqn:findxunderb} in step 6 replaced by
\begin{align*}
[\bar{x}_{k+1 \vert i} ]_j &= \! \max_{x\in X^0_{k|i}} 
[g(x,Kx \!+\! c_{k|i}^0)]_{j} \! - [h(x_{k|i, j}^{0,h},u_{k|i,j}^{0,h})]_{j} +[\bar{w}]_j 
\\
[\munderbar{x}_{k+1 \vert i} ]_j &= 
[g(x_{k|i, j}^{0,g},u_{k|i,j}^{0,g})]_{j} - \!\!
\min_{x\in X^0_{k|i}} [h(x,Kx \!+\! c_{k|i}^0)]_{j}  + [\munderbar{w}]_j 
\end{align*}     
and with constraints \eqref{eqn:RcMPCRHOCP4} in the RcMPC RHOCP replaced by 
\begin{align}
[\bar{x}_{k+1 \vert i}]_j &\geq [g(x, Kx+c_{k|i} + c^0_{k|i})]_j - [h(x^{0,h}_{k \vert i,j},u^{0,h}_{k \vert i,j})]_{j} \notag \\
& \quad - [A^h_{k \vert i} \!+\! B^h_{k \vert i}K]_j (x\!-x^{0,h}_{k \vert i,j}) - [B^h_{k \vert i}]_j c_{k \vert i} + [\bar{w}]_j
\nonumber\\
[\munderbar{x}_{k+1 \vert i}]_j &\leq [g(x^{0,g}_{k \vert i,j},u^{0,g}_{k \vert i,j})]_{j} + [A^g_{k \vert i} \!+\! B^g_{k \vert i}K]_j (x\! -x^{0,g}_{k \vert i,j}) \notag  \\ 
&\quad + [B^g_{k \vert i}]_j c_{k \vert i}  - [h(x, Kx+c_{k|i} + c^0_{k|i})]_j
+ [\munderbar{w}]_j.
\nonumber
\end{align}
\end{rem}

  \begin{rem}\label{rem:simplex_tube}
    The RcMPC RHOCP can be implemented by replacing $x$ in \eqref{eqn:RcMPCRHOCP4} by each vertex of $X_{k|i}$. If $X_{k|i}$ is defined by elementwise bounds, as in \eqref{eq:Xsetdef}, the number of constraints  will grow exponentially with the state dimension $n_x$. However, this exponential growth can be avoided if ${\bf X}_i$ is defined instead as a homothetic tube.
    \textcolor{black}{For example, any closed convex polytopic tube cross-section can be expressed  
    \begin{equation}
        X_{k \vert i}:= \{x: Qx\leq q_{k \vert i}\} = {\mathrm co}\{Z^v q_{k \vert i}, v = 1,\ldots,n_v\}, \notag
    \end{equation}
    where $q_{k \vert i}$ is a variable, while $Q,Z^1,\ldots,Z^{n_v}$ are fixed and ${\mathrm co}$ denotes the convex hull. 
    Writing $Qf(x,u) = \tilde{g}(x,u) - \tilde{h}(x,u)$, where $\tilde{g},\tilde{h}$ are differentiable and componentwise convex, we define $x^{0,\tilde{h}}_{k|i,j}$ as in Remark~\ref{rem:DC_dist} and replace~\eqref{eqn:RcMPCRHOCP4} by
    \begin{multline}\label{eqn:two22}
    [q_{k+1 \vert i}]_j \geq [\tilde{g}(x, Kx+c_{k|i} + c^0_{k|i})]_j - [\tilde{h}(x^{0,\tilde{h}}_{k \vert i,j},u^{0,\tilde{h}}_{k \vert i,j})]_{j} \\
    -\! [A^{\tilde{h}}_{k \vert i} \!+\! B^{\tilde{h}}_{k \vert i}K]_j (x\!-x^{0,\tilde{h}}_{k \vert i,j}) \!-\! [B^{\tilde{h}}_{k \vert i}]_j c_{k \vert i} \!+\! \max_{w \in \mathcal{W}} [Q]_j w
    \end{multline}
    with $x=Z^v q_{k|i}$, $v= 1,\ldots,n_v$.
%
This is identical to~\eqref{eqn:RcMPCRHOCP4} if $Q = [I \ {-I}]^\intercal$, $q_{k|i} = [\bar{x}_{k|i}^\intercal \ {-\munderbar{x}_{k|i}^\intercal}]^\intercal$, 
but if $Q=\bigl[-I \ [1 \cdots 1]^{\!\intercal}\bigr]^{
\!\intercal}$, we have simplex cross-sections with only $n_v=n_x+1$ vertices.}%
\end{rem}%

\subsection{Robust feasibility and stability guarantees }
\label{section:RobustFeasStabGuar} 

This section shows that, given an initial perturbation sequence $\mathbf{c}_0^0$ such that Algorithm \ref{alg:cap4} with $x_0=x_{\mathrm{init}}$ generates a feasible tube, the RcMPC RHOCP remains feasible at each iteration of Algorithm~\ref{alg:cap2} at all times. 
We also show that the optimal objective of the RcMPC RHOCP is non-increasing, and therefore the iteration of Algorithm~\ref{alg:cap2} converges and the MPC law provides a robust quadratic stability guarantee for any chosen limit on the number of iterations satisfying $\mathit{maxiters}\geq 1$.

To make notation more precise, at the $i$th time step and $n$th iteration of Algorithm~\ref{alg:cap2} we denote 
the seed trajectory as $\smash{\mathbf{c}_i^{0,n}}$ and the tube generated by Algorithm~\ref{alg:cap4} with $\smash{\mathbf{c}_i^{0}}={\mathbf{c}_i^{0,n}}$ as $\smash{\mathbf{X}_i^{0,n}}$. 
Similarly, we denote the solution of the RcMPC RHOCP in line~4  as $(\smash{\mathbf{c}_i^{\ast,n}, \mathbf{X}_i^{\ast,n}})$, and the optimal objective as $\smash{J_{\mathrm{RcMPC},i}^{\ast,n}} = \smash{J_{\mathrm{RcMPC}}(\mathbf{c}_i^{\ast,n}, \mathbf{X}_i^{\ast,n}, \mathbf{c}_i^{0,n})}$.

First note that the tubes ${\mathbf{X}_i^{0,n}}$ and ${\mathbf{X}_i^{\ast,n}}$ generated by Algorithms~\ref{alg:cap4} and \ref{alg:cap2} necessarily contain all trajectories of \eqref{eqn:robustsystemmodel} with initial condition $x_{0|i}=x_i$ under the control laws $u_{k|i} = Kx_{k|i} + \smash{c_{k|i}^{0,n}}$ and $u_{k|i} = Kx_{k|i} + \smash{c_{k|i}^{\ast,n}}$ respectively. This follows from the convexity of $f$ (Assumption~\ref{assumpt:two}), and from the definition of $\smash{X_{k+1|i}^{0,n}}$ in \eqref{eqn:findxoverb}-\eqref{eqn:findxunderb}, and the constraints on $\smash{X_{i+1|k}^{\ast,n}}$ in \eqref{eqn:RcMPCRHOCP4}, which imply
\begin{alignat}{2}%
& f(x,Kx+c_{k|i}^{0,n})+w \in X_{k+1|i}^{0,n} &\quad & \forall x \in X_{k|i}^{0,n}
\\
& f(x,Kx+c_{k|i}^{\ast,n})+w \in X_{k+1|i}^{\ast,n} &\quad & \forall x \in X_{k|i}^{\ast,n}
\end{alignat}%
for all $w\in\W$, $k=0,\ldots,N-1$ and $i\geq 0$, where 
$\smash{\mathbf{c}_i^{0,n}} = \{\smash{c_{0|i}^{0,n}},\ldots,\smash{c_{N-1|i}^{0,n}}\}$ and 
$\smash{\mathbf{c}_i^{\ast,n}} = \{\smash{c_{0|i}^{\ast,n}},\ldots,\smash{c_{N-1|i}^{\ast,n}}\}$. 

The recursive feasibility of Algorithm~\ref{alg:cap2} is a consequence of the following properties of the sets $X_{i|k}^{0,n}$, $X_{i|k}^{\ast,n}$.

\begin{lem}\label{lem:tubewithintube1}%
For all $i\geq 0$ and $n=1,\ldots n_i$ we have
\begin{alignat}{2}
&\text{if } X_{k|i}^{0,n+1} \subseteq X_{k|i}^{\ast,n}
&\quad & 
\text{then } X_{k+1|i}^{0,n+1}\subseteq X_{k+1|i}^{\ast,n}
\label{eqn:tint1}\\
&\text{if } X_{k|i+1}^{0,1}\subseteq X_{k+1|i}^{\ast,n_i}
&\quad & 
\text{then } X_{k+1|i+1}^{0,1}\subseteq X_{k+2|i}^{\ast,n_i}
\label{eqn:tint2}\\
&\text{if } X_{N-1|i+1}^{0,1}\subseteq X_{N|i}^{\ast,n_i}
&\quad & 
\text{then } X_{N|i+1}^{0,1}\subseteq \X_N
\label{eqn:tint3}
\end{alignat}
where \eqref{eqn:tint1} holds for $k=0,\ldots,N-1$, \eqref{eqn:tint2} holds for $k=0,\ldots,N-2$, and $n_i$ is the final iteration at time $i$.
\end{lem}

\textbf{Proof.}\hspace{1ex}%
Condition \eqref{eqn:tint1} follows from the definition of the linearization points in \eqref{eqn:findlinearizpoints}, which implies
\begin{alignat}{1}
&
\text{\makebox[0pt][l]{$[f(x_{k|i,j}^{0,n+1},u_{k|i,j}^{0,n+1})]_j
= \min_{x\in X_{k|i}^{0,n+1}} [f(x,Kx+c_{k|i}^{0,n+1})]_j$}}
\label{eqn:tintsub1}\\
&\quad \geq \min_{x\in X_{k|i}^{\ast,n}} [f(x,Kx+c_{k|i}^{0,n+1})]_j
\label{eqn:tintsub2}\\
&\quad\geq \min_{x\in X_{k|i}^{\ast,n}} \Bigl\{ [f(x_{k|i,j}^{0,n},u_{k|i,j}^{0,n}]_j + [A_{k|i}^{n}]_j(x-x_{k|i,j}^{0,n}) 
\nonumber\\
&\qquad\qquad\quad + [B_{k|i}^{n}]_j (Kx + c_{k|i}^{0,n+1} - u_{k|i}^{0,n}) \Bigr\}
\label{eqn:tintsub3}\\
&\quad= 
\min_{x\in X_{k|i}^{\ast,n}} \Bigl\{ [f(x_{k|i,j}^{0,n},u_{k|i,j}^{0,n}]_j \!+\! [A_{k|i}^{n}\!+\!B_{k|i}^{n}K]_j(x\!-\!x_{k|i,j}^{0,n}) 
\nonumber\\
&\qquad\qquad\quad + [B_{k|i}^{n}]_j c_{k|i}^{\ast,n}\Bigr\}
\label{eqn:tintsub4}
\end{alignat}
where \eqref{eqn:tintsub2} holds if $\smash{X_{k|i}^{0,n+1}}\subseteq \smash{X_{k|i}^{\ast,n}}$, \eqref{eqn:tintsub3} is due to the convexity of $f$, and \eqref{eqn:tintsub4} follows from $u_{k|i}^{0,n}=Kx_{k|i}^{0,n}+c_{k|i}^{0,n}$ and $\smash{c_{k|i}^{0,n+1}} = \smash{c_{k|i}^{\ast,n}}+\smash{c_{k|i}^{0,n}}$. Therefore, from \eqref{eqn:RcMPCRHOCP4} we obtain
\[
\munderbar{x}_{k+1|i}^{0,n+1} \geq \munderbar{x}_{k+1|i}^{\ast,n}.
\]
Similarly, if $X_{k|i}^{0,n+1}\subseteq X_{k|i}^{\ast,n}$, then
\begin{gather*}
\max_{x\in X_{k|i}^{0,n+1}} [f(x, Kx+c_{k|i}^{0,n+1}]_j 
\leq 
\max_{x\in X_{k|i}^{\ast,n}} [f(x, Kx+c_{k|i}^{0,n+1}]_j 
\\
=
\max_{x\in X_{k|i}^{\ast,n}} [f(x, Kx+c_{k|i}^{\ast,n}+c_{k|i}^{0,n}]_j 
\end{gather*}
and hence
\[
\bar{x}_{k+1|i}^{0,n+1} \leq \bar{x}_{k+1|i}^{\ast,n}
\]
by \eqref{eqn:findxoverb} and \eqref{eqn:RcMPCRHOCP4}, which demonstrates \eqref{eqn:tint1}.
Conditions \eqref{eqn:tint2} and \eqref{eqn:tint3} follow from the $\smash{\mathbf{c}_{i+1}^{0}}$ update in line 6 of Algorithm~\ref{alg:cap2}, which implies $\smash[b]{c_{k|i+1}^{0,1}} = \smash[b]{c_{k+1|i}^{\ast,n_i}}+\smash[b]{c_{k+1|i}^{0,n_i}}$, $k=0,\ldots,N-2$, and $\smash[b]{c_{N-1|i+1}^{0,1}} = u_{N+i}^r-Kx_{N+i}^r$. Hence \eqref{eqn:tint2} can be shown analogously to \eqref{eqn:tint1}, and \eqref{eqn:tint3} holds because $X_{N|i}^{\ast,n_i}\subseteq \X_N$ by \eqref{eqn:RcMPCRHOCP5} and $\X_N$ satisfies Assumption~\ref{assumpt:robust}. 
\qed

\begin{thm}[Recursive feasibility] \label{thm:robustrecurfeas}
If \mbox{Assumptions} \ref{assumpt:two} and \ref{assumpt:robust} hold, 
and if $\mathbf{c}_0,\mathbf{X}_0$ exist satisfying \eqref{eqn:RcMPCRHOCP3}-\eqref{eqn:RcMPCRHOCP5} with the initial seed $\mathbf{c}_0^0$, then the RcMPC RHOCP is feasible at each iteration of Algorithm~\ref{alg:cap2}, for all times $i\geq 0$.
\end{thm}

\textbf{Proof.}\hspace{1ex}%
If, at time $i$ and iteration $n$ the RcMPC RHOCP in line~4 of Algorithm~\ref{alg:cap2} has solution $(\smash{\mathbf{c}_i^{\ast,n}},\smash{\mathbf{X}_i^{\ast,n}})$, then $(\smash{\mathbf{c}_i},\smash{\mathbf{X}_i})= (0,\smash{\mathbf{X}_i^{0,n+1}})$ is a feasible solution at iteration ${n+1}$ by Lemma~\ref{lem:tubewithintube1} since $\smash{X_{0|i}^{0,n+1}}=\{x_i\}=\smash{X_{0|i}^{\ast,n}}$.
Similarly, Lemma~\ref{lem:tubewithintube1} implies that, if $(\smash{\mathbf{c}_i^{\ast,n_i}},\smash{\mathbf{X}_i^{\ast,n_i}})$ is the solution of the RcMPC RHOCP at the final iteration performed at time $i$, then $(\smash{\mathbf{c}_{i+1}},\smash{\mathbf{X}_{i+1}})= (0,\smash{\mathbf{X}_{i+1}^{0,1}})$ is a feasible solution at the first iteration at time $i+1$ since ${X_{0|i+1}^{0,1}}=\{x_{i+1}\}\subseteq {X_{0|i}^{\ast,n_i}}$.
Furthermore, feasibility of \eqref{eqn:RcMPCRHOCP3}-\eqref{eqn:RcMPCRHOCP5} with the initial seed $\smash{\mathbf{c}_0^0}$ ensures feasibility of the RcMPC RHOCP at the first iteration at $i=0$, so feasibility is ensured at all iterations and all time steps.
\qed

\begin{lem} \label{lem:JRcmpc}
For all $i\geq 0$ and $n \geq 1$ the optimal cost satisfies $ J_{\mathrm{RcMPC},i}^{*,n+1} \leq J_{\mathrm{RcMPC},i}^{*,n}$.
\end{lem}

\textbf{Proof.}\hspace{1ex}
Lemma~\ref{lem:tubewithintube1} implies $(\smash{\mathbf{c}_i},\smash{\mathbf{X}_i})= (0,\smash{\mathbf{X}_i^{0,n+1}})$ is a feasible solution of the RcMPC RHOCP at iteration $n+1$. Therefore the optimal solution at iteration ${n+1}$ satisfies
${J_{\mathrm{RcMPC},i}^{\ast,n+1}} \leq {J_{\mathrm{RcMPC}}(0, {\mathbf{X}_i^{0,n+1}},{\mathbf{c}_i^{0,n+1}})}$. But from Lemma~\ref{lem:tubewithintube1} we also have ${X_{k|i}^{0,n+1}}\subseteq {X_{k|i}^{\ast,n}}$ for $k=0,\ldots,N$, and the cost~\eqref{eqn:RcMPCRHOCP0} therefore ensures $J_{\mathrm{RcMPC}}(\textcolor{red}{0}, {\mathbf{X}_i^{0,n+1}},\textcolor{red}{{\mathbf{c}_i^{\ast,n}}})\leq J_{\mathrm{RcMPC},i}^{\ast,n}$.
\qed 

\begin{rem}\label{rem:feastubesolut}  
Since $J_{RcMPC,i}^{\ast,n}\geq 0$, Lemma~\ref{lem:JRcmpc} implies that, in the ideal case of unlimited iterations, the optimal RcMPC cost in line~4 of Algorithm~\ref{alg:cap2} converges to a finite limit as $n\to\infty$. As in the nominal case, it follows that $\smash{\mathbf{c}_i^{\ast,n}}\to 0$ as $n\to\infty$.
In this case however, the seed tube $\mathbf{X}_i^{0,n}$ generated by Algorithm~\ref{alg:cap4} will in general contain more than a single trajectory at all iterations $n\geq 1$, and hence $\smash{\mathbf{X}_i^{0,n}}$ does not converge to a single trajectory. 
\end{rem}

\begin{thm}\label{thm:robust_stability}
If Assumptions \ref{assumpt:two} and \ref{assumpt:robust} and hold and the offline computation of Algorithm~\ref{alg:cap2} is feasible, then the closed loop system formed by \eqref{eqn:robustsystemmodel} and 
Algorithm~\ref{alg:cap2} 
satisfies the state and input constraints \eqref{eqn:stateandinputconstr} and the quadratic stability condition: 
\begin{equation}\label{eqn:limitcond}
\lim_{t \to \infty} \frac{1}{t} \sum_{i=0}^{t-1} ( \| x_{i} - x_{i}^r \|_Q^2 + \| u_{i} - u_{i}^r \|_R^2 ) \leq \beta 
\end{equation}
where $\beta$ satisfies Assumption~\ref{assumpt:robust}(c).  
\end{thm}

\textbf{Proof.}\hspace{1ex}
By Lemma~\ref{lem:tubewithintube1} and Theorem~\ref{thm:robustrecurfeas}, $(\mathbf{c}_{i+1},\mathbf{X}_{i+1})=(0,\mathbf{X}_{i+1}^{0,1})$ is a feasible solution for the RcMPC RHOCP at the first iteration of Algorithm~\ref{alg:cap2} at time $i+1$, so
\begin{align*}
J_{\mathrm{RcMPC},i+1}^{\ast,1} 
&\leq 
J_{\mathrm{RcMPC}}(0,\mathbf{X}_{i+1}^{0,1},\mathbf{c}_{i+1}^{0,1}) \\
&\leq 
J_{\mathrm{RcMPC}}(\mathbf{c}_i^{\ast,n_i},\mathbf{X}_{i}^{\ast,n_i},\mathbf{c}_{i}^{0,n_i}) 
\\
&\quad
- \|x_i- x_i^r\|_Q^2 - \|u_i - u_i^r\|_i^2 + \beta,
\end{align*}
where the definition of $\mathbf{c}_{i+1}^{0,1}$ in line 6 and the inequality in Assumption~\ref{assumpt:robust}(c) have been used.
But Lemma~\ref{lem:JRcmpc} implies
$\smash{J_{\mathrm{RcMPC},i+1}^{\ast,n_{i+1}}}
\leq 
\smash{J_{\mathrm{RcMPC},i+1}^{\ast,1}}$ and hence
\begin{align}\label{eq:iss_lyap}
J_{\mathrm{RcMPC},i+1}^{\ast,n_{i+1}} &\leq  
J_{\mathrm{RcMPC},i}^{\ast,n_{i}} - \|x_i- x_i^r\|_Q^2 - \|u_i - u_i^r\|_i^2 + \beta .
\end{align}
The bound \eqref{eqn:limitcond} follows by summing both sides of \eqref{eq:iss_lyap}
over $i= 0,\ldots,t-1$ and 
noting that $\smash[t]{J_{\mathrm{RcMPC},i}^{\ast,n_i}}$ is  finite for all $i$ since $\X$, $\U$ are compact and the RcMPC RHOCP is feasible.
\qed 

\begin{rem}
If Assumption~\ref{assumpt:robust} reduces to Assumption~\ref{assumpt:matrices} in the absence of disturbances (i.e.~if $\beta\to 0$ as $\mathcal{W}\to \{0\}$), then the origin ($x-x^r = 0$) of the closed loop system of \eqref{eqn:robustsystemmodel} with Alg.~\ref{alg:cap2} 
is input-to-state stable (ISS) \cite{limon2009iss}, since \eqref{eq:iss_lyap} and Assumption~\ref{assumpt:robust}(c) imply 
$\smash{J_{\mathrm{RcMPC},i}^{\ast,n_i}}$ is an ISS Lyapunov function.
\end{rem}

\section{Numerical examples}
This section presents numerical results obtained by applying the nominal cMPC (Algorithm \ref{alg:cap}) and robust RcMPC (Algorithm \ref{alg:cap2}) strategies to two example problems. In each case the algorithms are implemented in Python, using cvxpy \cite{cvxpy} and MOSEK~\cite{mosek} for convex optimisation problems, and using the SciPy SLSQP algorithm~\cite{scipy} for comparisons with a general nonlinear programming solver. 
Unless stated otherwise $\textit{maxiters} = 3$, $\textit{tol} = 10^{-6}$ and MPC simulations are terminated when $\| x_{N \vert i} \|_2$ reaches $10^{-4}$ (Example 1) or $10^{-6}$ (Example~2).

\textbf{Example 1}\hspace{1ex} First we consider a system with model
\begin{align*}
    \dot{x}_{1}(t)&= x_2(t) +w_1(t) \\
    \dot{x}_{2}(t) &= 0.2e^{-x_1(t)}-x_2(t)+u_1(t)-0.2+w_2(t) .
\end{align*}
and cost weights $Q = \textrm{diag}\{1,1\} $, $R = 1$, and initial condition $x_1(0) = 5$, $x_2(0) = 10$. Forward Euler discretization is used with time step $\Delta t = 8\times10^{-3}$, and the prediction horizon is $N=25$. The reference trajectory is chosen as $x^r_i = 0, u^r_i = 0$ for all $i$, and the system constraints are $\vert x_1(t) \vert \leq 10$, $\vert x_2(t) \vert \leq 10$, $\vert u_1(t) \vert \leq 150$ at all discrete times $t=i\Delta t$, $i=0,1,\ldots$.

\begin{figure}[!tbp]
\centerline{\includegraphics[scale=0.16]{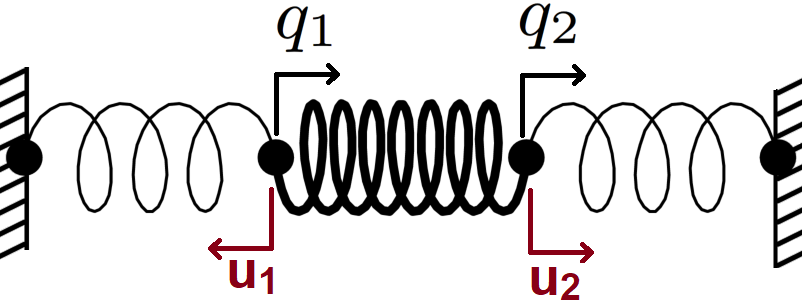}}
\caption{FPU system with two unit masses with coordinates $q_1,q_2$, and control forces $u_1,u_2$ applied directly to each mass.}
\label{fig:SPRINGS}
\end{figure}

\textbf{Example 2}\hspace{1ex} Applying a change of coordinates to the Fermi Pasta Ulam (FPU) chain with cubic spring potentials (Figure \ref{fig:SPRINGS}), we obtain a convex model formulation:
 \begin{align*}
    \dot{x}_{1}(t) &= x_3(t) +w_1(t) \\
    \dot{x}_{2}(t) &= x_4(t) +w_2(t) \\
    \dot{x}_{3}(t) &= -\tfrac{\eta^2}{2}\Big(x_2(t)+x_1(t)\Big)+3x_1(t)^2 + u_1(t)+ w_3(t) \\
    \dot{x}_{4}(t) &= -\tfrac{\eta^2}{2}\Big(x_2(t)+x_1(t)\Big)+3x_2(t)^2 +u_2(t) + w_4(t) .
\end{align*}
Forward Euler discretization with $\Delta t = 10^{-2}$ and a prediction horizon $N=25$ are used. The model and algorithm parameters are $\eta = 50$,  $x^r_i = 0$, $u^r_i = 0$ for all $i$, $Q = 10^{-2} I$, $R = 10^{-2} I$, 
$\lvert x_1(t) \rvert \leq 10$, $\lvert x_2(t) \rvert \leq 10$, $\lvert x_3(t) \rvert \leq 10$, $\lvert x_4(t) \rvert \leq 10$, $\lvert u_1(t) \rvert \leq 33$, $\lvert u_2(t) \rvert \leq 33$, $t=i\Delta t$, $i=0,1,\ldots$, and
\[
\begin{bmatrix}
   x_1(0) \\ x_2(0)
\end{bmatrix} = E \begin{bmatrix}
    1 \\ 1
\end{bmatrix} \! , 
\begin{bmatrix}
   x_3(0) \\x_4(0)
\end{bmatrix} = E \begin{bmatrix}
    1/\eta \\ 1
\end{bmatrix} \! , 
E   =\tfrac{1}{\sqrt{2}}\begin{bmatrix}
   - 1 & 1 \\1&1  
\end{bmatrix}.
\]

\subsection{Nominal cMPC}

\begin{figure}[!tbp] \centerline{\includegraphics[width=0.5\columnwidth]{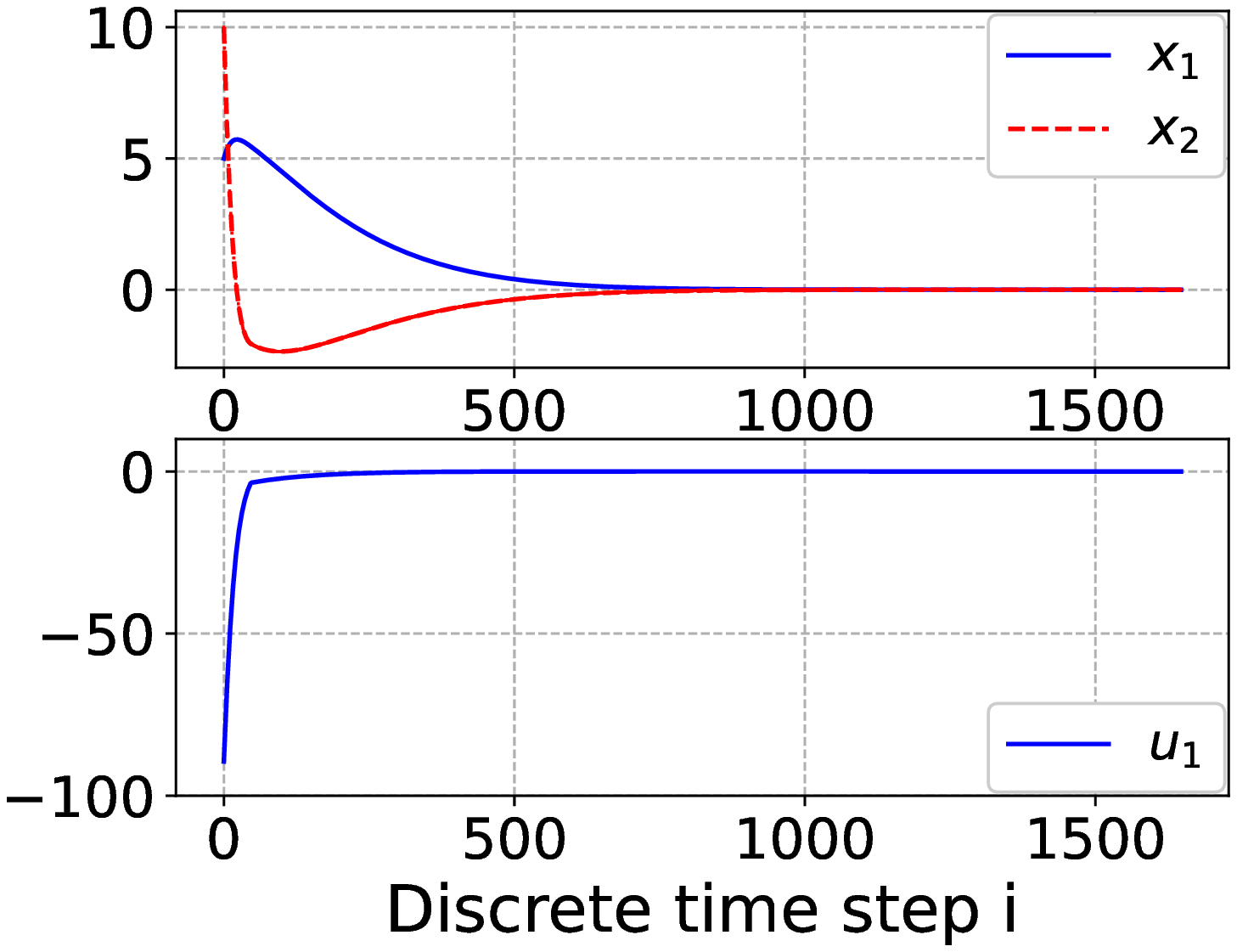}%
\includegraphics[width=0.49\columnwidth]{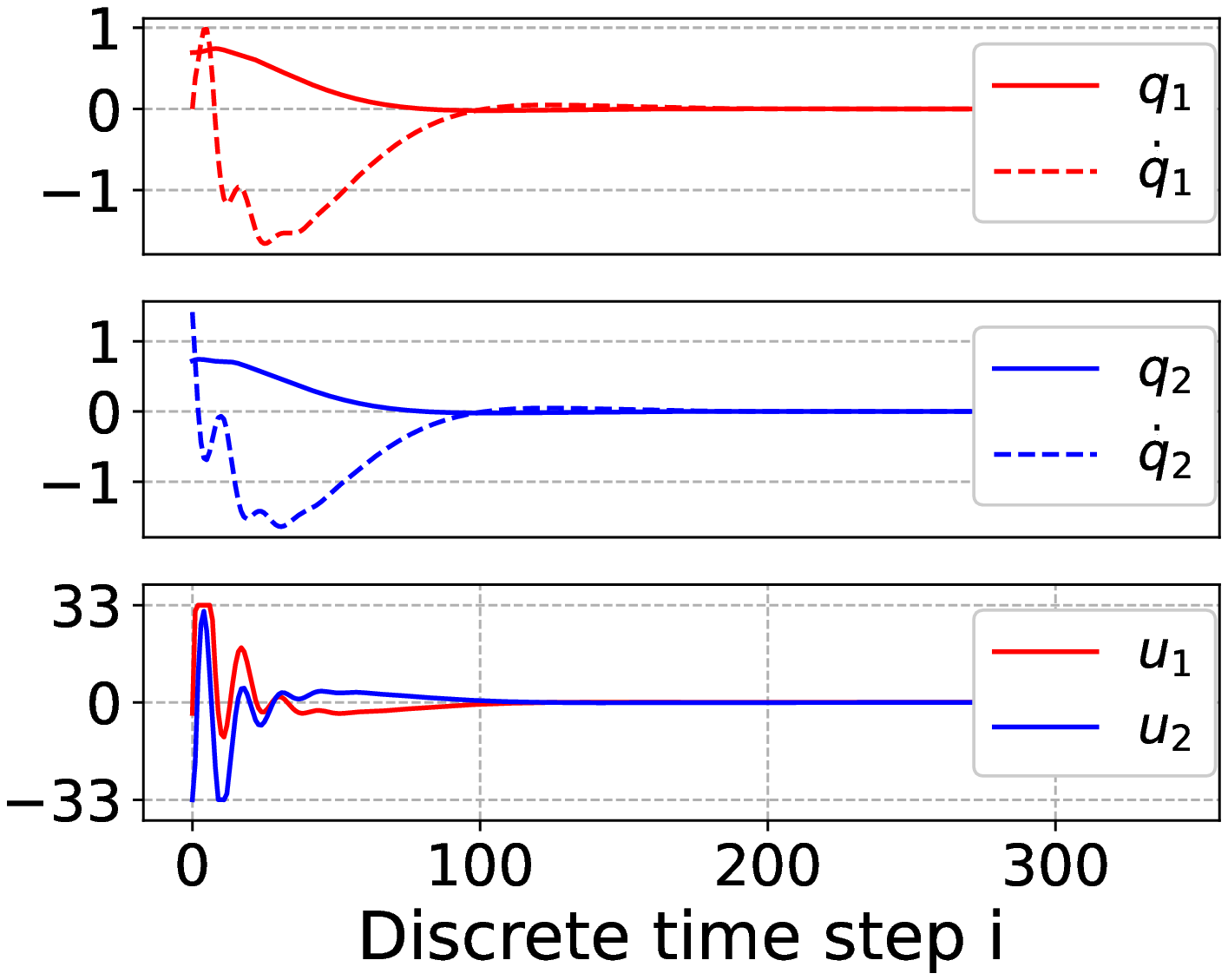}}
\caption{Closed loop trajectories with cMPC (Algorithm~\ref{alg:cap}). Left: Example 1. Right: Example 2.}
\label{fig:trajplot}
\end{figure}

\begin{figure}[!tbp]
\centerline{\includegraphics[width=0.9\columnwidth]{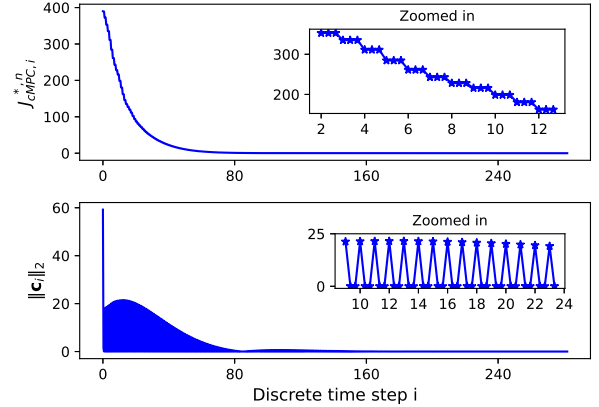}}
\caption{Convergence of cMPC at successive iterations and time steps for Example 2. Upper plot: Optimal cost. Lower: $\|\mathbf{c}_i^\ast\|$.}
\label{fig:controlplot}
\end{figure}

\begin{figure}[!tbp]
\centerline{\includegraphics[trim={3mm 3mm 1mm 0}, clip,width=1.05\columnwidth]{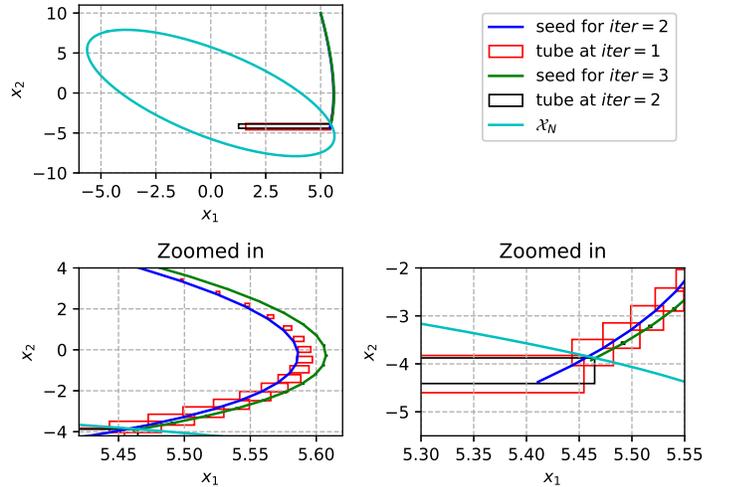}}
\caption{Seed trajectories and tubes (${x}_{i|k}^0$ and $x_{i|k}^0 \oplus S_{i|k}^\ast$) for cMPC and Example 1 at iterations $n=1,2,3$ at time $i=0$.}
\label{fig:tubeconvrg}
\end{figure}

The performance of the nominal cMPC strategy (Algorithm~\ref{alg:cap}) described in Section \ref{section:nominalsection} was investigated by setting $w_1(t)=w_2(t)=0$, $w_3(t)=w_4(t)=0$ for all $t$ in Examples 1 and 2. The initial seed was obtained using the initialization strategy of Section \ref{section:nominalsection}. 
Figure~\ref{fig:trajplot} shows the predicted trajectories of Algorithm \ref{alg:cap} for both example problems. 
As expected from Lemma~\ref{lem:cost_sequence} and Theorem~\ref{thm:convergence}, the optimal predicted cost is non-increasing at successive iterations and time steps, and $\mathbf{c}_i^\ast$ converges to $0$. This is shown in Fig.~\ref{fig:controlplot} for Example~2, where the points between integer time steps represent the iterations performed at each time step.
The convergence of the tube to a single trajectory (Theorem~\ref{thm:convergence} and Remark~\ref{rem:convergence}) is illustrated in Fig.~\ref{fig:tubeconvrg} for Example~1. The iteration  converges rapidly ($\mathbf{S}_0^\ast$ is already difficult to see at the second iteration) and the tube at iteration $n$ contains the seed trajectory updated for iteration $n+1$ via line 6 of Algorithm~\ref{alg:cap}, as expected from Theorem~\ref{thm:recrsfeas}.

\begin{table}[h]
\resizebox{\columnwidth}{!}{%
\begin{tabular}{@{\hspace{0cm}}l@{\hspace{0.05cm}}l@{\hspace{-0.15cm}}c@{\hspace{-0cm}}c@{\hspace{0.25cm}}c@{\hspace{0.25cm}}c@{\hspace{0.25cm}}c@{\hspace{0.25cm}}c@{\hspace{0.25cm}}c@{\hspace{0.25cm}}c@{\hspace{0.25cm}}c@{\hspace{0.25cm}}c@{\hspace{0cm}}}
\hline
& \multirow{4}{*}{\hspace{0.1cm}\rot{60}{solver}} & \multicolumn{1}{l}{\multirow{4}{*}{\rot{60}{maxiters}}} 
&  &  &  & time & \textbf{total} & total & total & \multirow{4}{*}{$J_0^\ast$} & \multirow{4}{*}{$J_{cl}$} \\
 &  & \multicolumn{1}{l}{} & $n_v$ & $n_c$ & $n_p$ & \multicolumn{1}{l}{for} & \textbf{time} & time & time &  &  \\
 &  & \multicolumn{1}{l}{} &  &  &  & \multicolumn{1}{l}{iter 1} & \textbf{at} & at & at &  &  \\
 &  & \multicolumn{1}{l}{} &  &  &  & \multicolumn{1}{l}{$i=0$} & \textbf{$i=0$} & $i=1$ & $i=2$ &  &  \\ \hline
\multirow{8}{*}{\rot{90}{Example 1}} & cvxpy& 1 & \multirow{6}{*}{180} & \multirow{6}{*}{510} & \multirow{6}{*}{277} & 0.038 & 0.229 & 0.053 & 0.055 & 132297 & 80588 \\
 & cvxpy & 3 &  &  &  & 0.039 & 0.331 & 0.168 & 0.154 & 121999 & 80152 \\
 & cvxpy & 5 &  &  &  & 0.040 & 0.494 & 0.258 & 0.262 & 121932 & 80138  \\
 & cvxpy$+$jit & 1 &  &  &  & 0.040 & \textbf{0.255} & \textbf{0.039} & \textbf{0.057} & 132297 &  80588 \\
 & cvxpy$+$jit & 3 &  &  &  & 0.039 & 0.292 & 0.124 & 0.122 & 121999 & 80152 \\
 & cvxpy$+$jit & 5 &  &  &  & 0.041 & \textbf{0.391} & \textbf{0.204} & \textbf{0.191} & 121932 &80137  \\ \cline{2-12} 
 & NLP & - & \multirow{2}{*}{77} & \multirow{2}{*}{355} & \multirow{2}{*}{-} & - & 13.401 &  17.881 & 12.828 & 121782 &77340  \\
 & NLP$+$jit & - &  &  &  & - & \textbf{1.547} & \textbf{1.562} & \textbf{1.507} & 121782 &   77340\\ \hline
\multirow{6}{*}{\rot{90}{Example 2}} & cvxpy & 1 & \multirow{4}{*}{309} & \multirow{4}{*}{2034} & \multirow{4}{*}{854} &0.970  & 1.284 & 0.192 &0.192  & 386.77 & 198.85 \\
 & cvxpy & 3 &  &  &  & 0.902 &  2.056&  0.634&  0.629& 389.27 & 198.80\\
 & cvxpy$+$jit & 1 &   &  &  & 0.992 & \textbf{1.259} & \textbf{0.167} &  \textbf{0.152} & 389.77 & 198.85 \\
 & cvxpy$+$jit & 3 &  &  &  &0.974 & \textbf{1.608} & \textbf{0.447} & \textbf{0.459}  & 389.27 & 198.80 \\ \cline{2-12} 
 & NLP & - & \multirow{2}{*}{154} & \multirow{2}{*}{709} & \multirow{2}{*}{-} & - & 31.602 & 34.365 & 37.811 & 389.27 & 198.80  \\
 & NLP$+$jit & - &  &  &  & - & \textbf{2.375} & \textbf{2.099} &  \textbf{2.060}&  389.27 & 198.80 \\ \hline
\end{tabular}%
}
\vspace*{0.5pt}
\caption{Comparison of solution times and optimal costs.}
\label{tab:timings}
\end{table}

Timing tests were conducted by comparing the cvxpy implementation of Algorithm \ref{alg:cap} to a receding horizon strategy that solves the nMPC RHOCP at each time step using the scipy.optimize library and its SLSQP solver. The cMPC RHOCP in line 4 of Algorithm~\ref{alg:cap} is implemented using the \textit{parameters} functionality of the cvxpy library and is thus precompiled offline. The solver timings were tested both with and without pre-compilation using the jit functionality of the jax library~\cite{jax2018github}. For the cvxpy algorithm jit is used to speed up lines 3 and 6, whereas for scipy.optimize it is used to pre-compile the constraints, constraint Jacobians, cost function and cost gradient. Results obtained using
a i7-10700 CPU at 2.90GHz are shown in Table \ref{tab:timings}, where $n_v$, $n_c$, $n_p$ are the number of variables, constraints and parameters respectively, $J_0^\ast$ is the optimal cost at $i=0$ and $J_{cl}$ is the cost evaluated using closed loop system responses.  
For Example 1 with standard implementation at $i=0$  Algorithm~\ref{alg:cap} requires only 1.71\% and 3.69\%. of the nMPC solution time for $\textit{maxiters} = 1$ and $5$ respectively. Even when special just-in-time compilation with jit is used to speed up both cMPC and nMPC, Algorithm~\ref{alg:cap} still presents great improvement requiring 16.48\% and 25.27\% of the nMPC solution time for  $\textit{maxiters} = 1$ and $5$ respectively. This percentage reduces to 2.50\% and 13.06\% at $i=1$.
For Example 2 we have a larger state space and hence more constraints in the cMPC RHOCP. Despite this we see that for $i=0$  Algorithm~\ref{alg:cap} requires  53.01\% and 67.71\% of the NLP solution time with jit for both approaches and 4.06\% and 6.51\% without jit, for $\textit{maxiters} = 1$ and $3$ respectively. Again for $i>0$ these percentages reduce significantly further. Furthermore, as $\textit{maxiters}$ increases the optimal 
cMPC RHOCP cost decreases and converges to the optimal nMPC RHOCP cost, which is consistent with Theorem~\ref{thm:KKT}.

\begin{figure}[!tbp]
\centerline{\includegraphics[width=0.49\columnwidth]{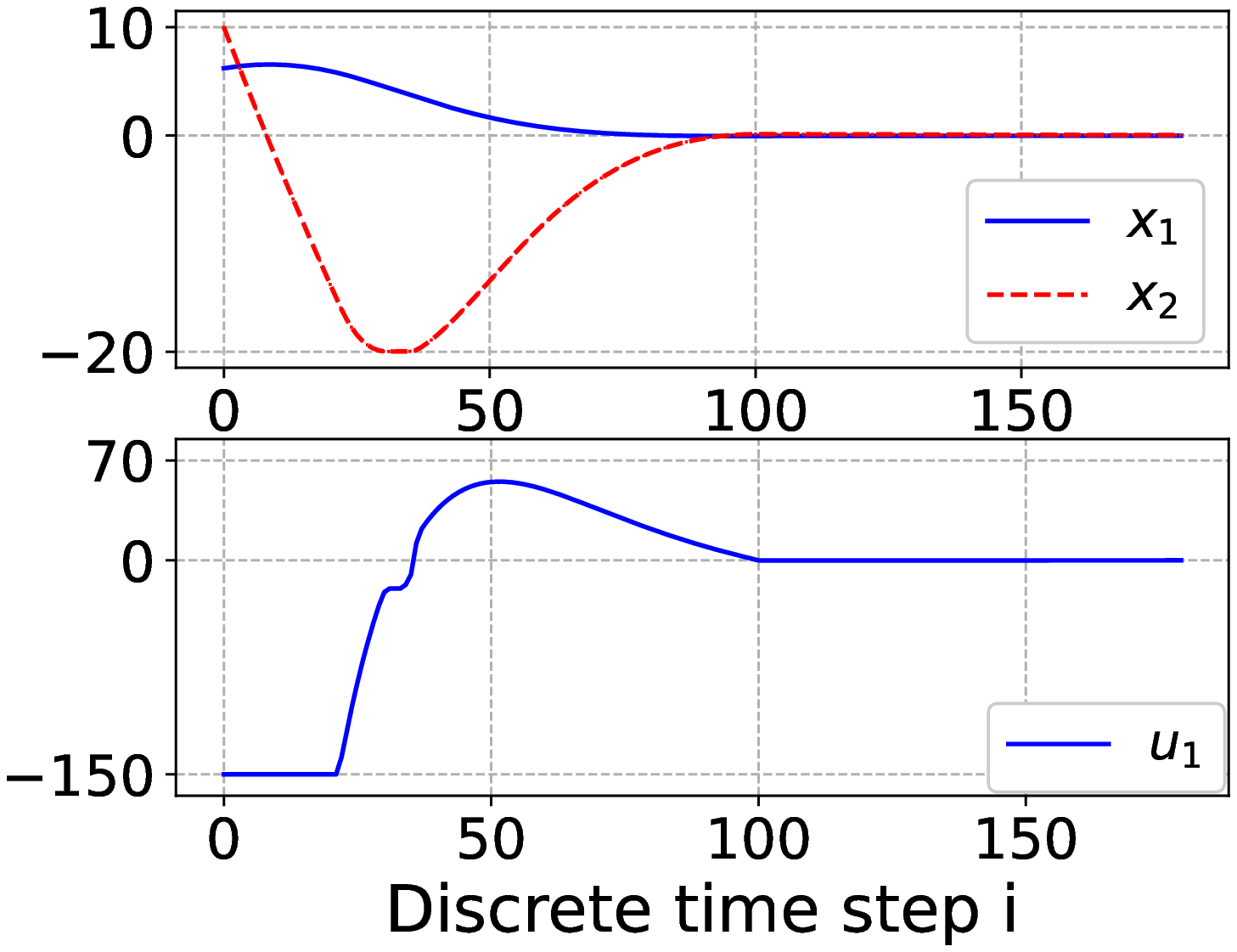}  
\includegraphics[width=0.49\columnwidth]{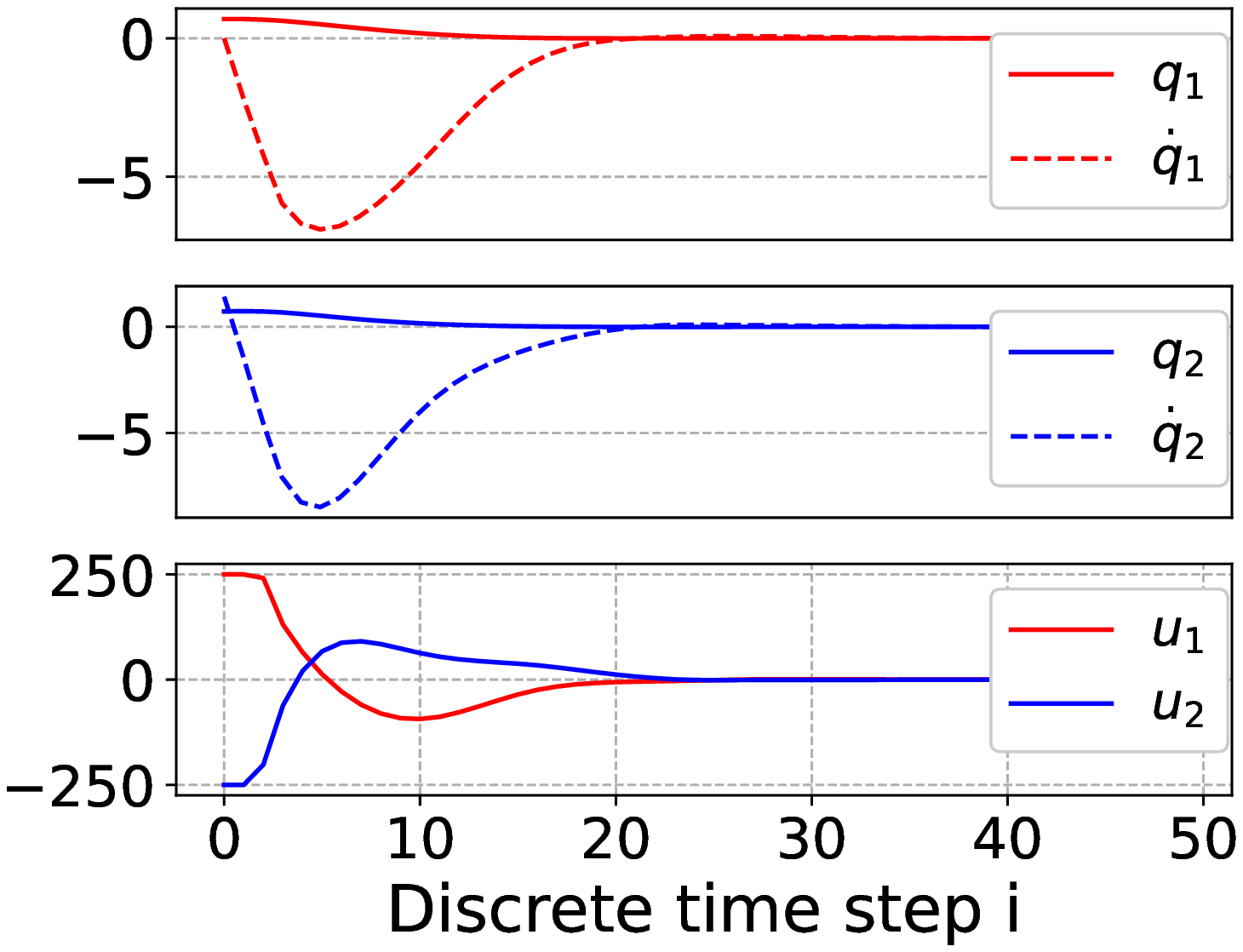}}
\vspace{-2mm}
\caption{Closed-loop state trajectories with RcMPC. Left: Example 1. Right: Example 2}
\label{fig:trajplot_R}
\end{figure}

\begin{figure}[!tbp]
\centerline{\includegraphics[trim={0 3mm 0 22mm},clip,width=0.49\columnwidth]{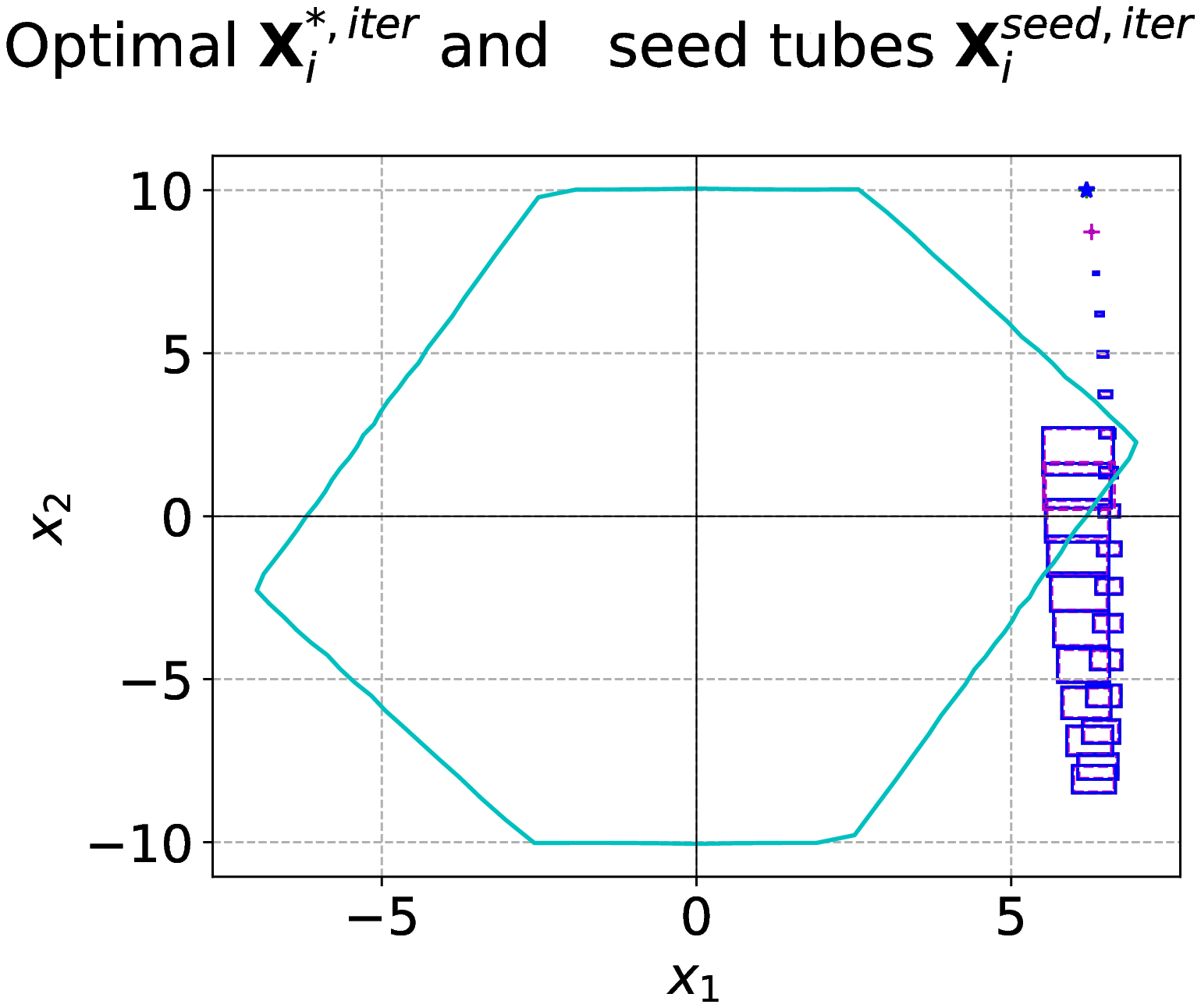}  
\includegraphics[trim={0 3mm 0 10mm},clip,width=0.49\columnwidth]{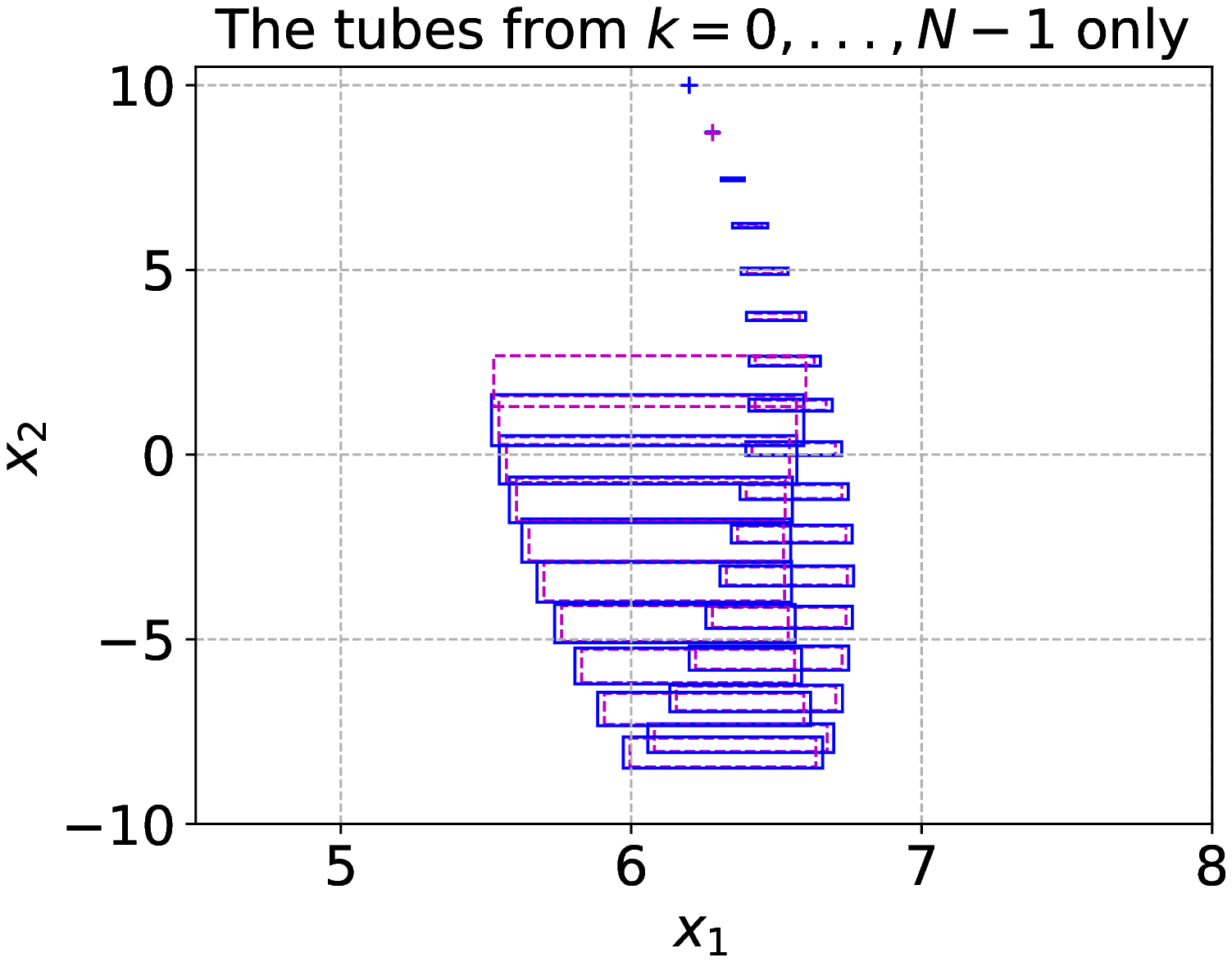}}
\vspace*{-2mm} 
\caption{Optimal and seed tubes for RcMPC for Example 1: $\smash{\mathbf{X}_0^{\ast,3}}$ (blue, solid lines), $\smash{\mathbf{X}_1^{0,1}}$ (pink dashed). Left: Tube cross-sections for prediction times $k=0,\ldots,N-1$. Right: Tube cross-sections for $k=0,\ldots,N$ and terminal set $\X_N$ (cyan).}
\label{fig:RcMPCplot1}
\end{figure}

\begin{figure}[!tbp]
\centerline{\includegraphics[width=0.5\columnwidth]{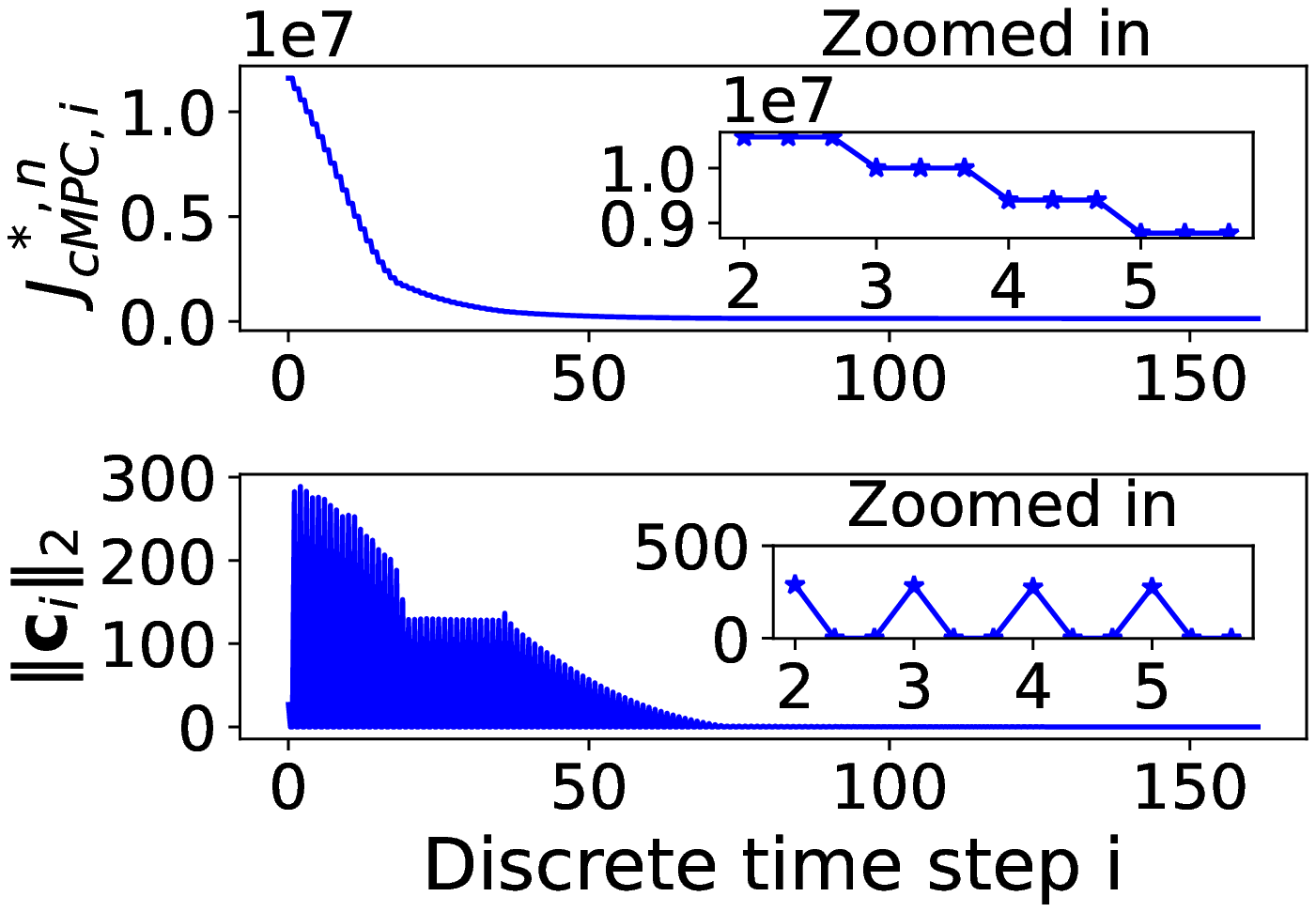}  
\includegraphics[width=0.5\columnwidth]{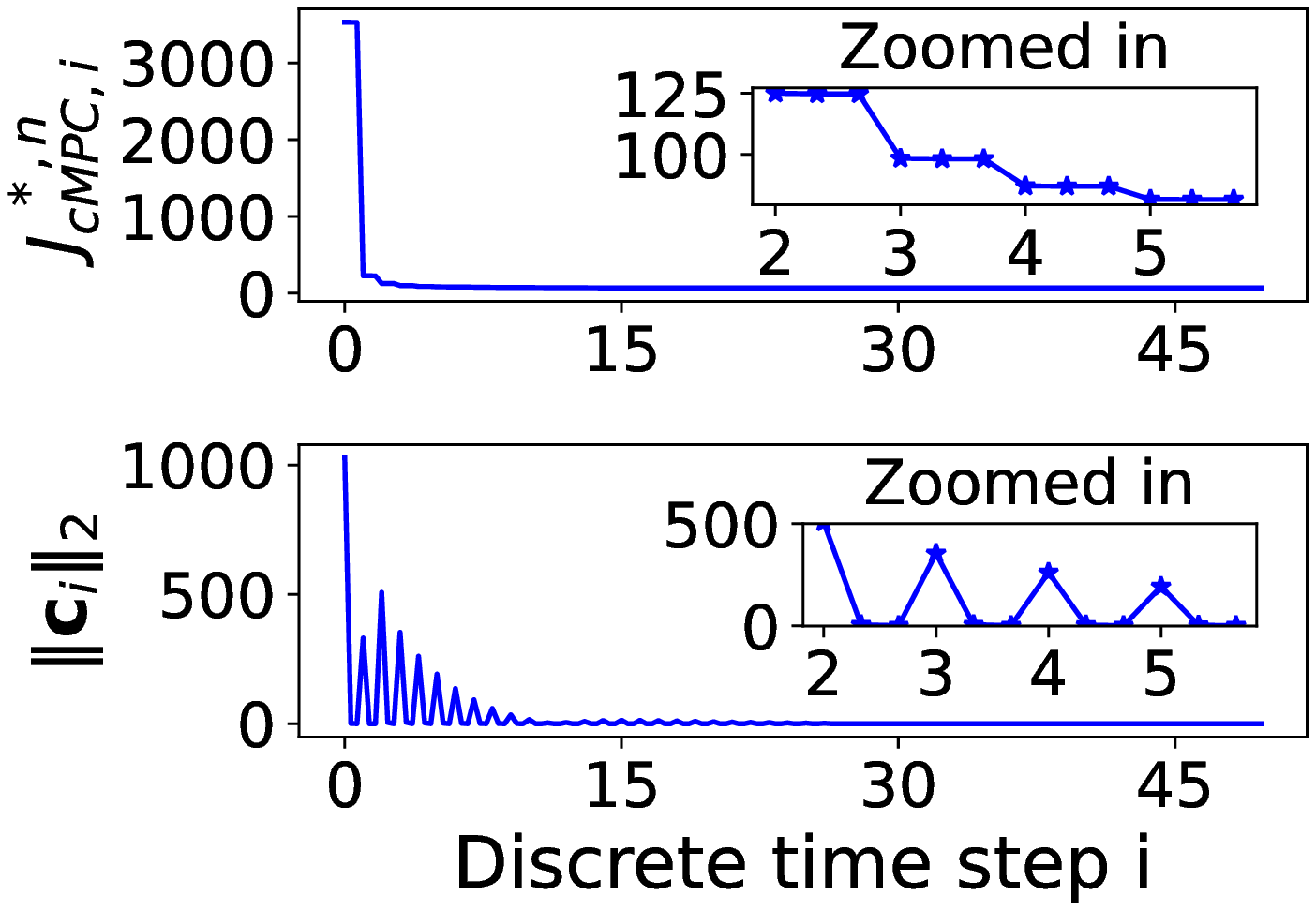}}
\vspace{-2mm}
\caption{Convergence of RcMPC optimal cost and control perturbation sequence $\|\mathbf{c}_i^\ast\|$. Left: Example 1. Right: Example~2.}
\label{fig:costplot_R1}
\vskip-\baselineskip
\end{figure}

\subsection{Robust cMPC (RcMPC)}
The performance of RcMPC (Algorithm~\ref{alg:cap2}) presented in Section \ref{section:robustsection} was investigated using Examples~1 and~2 with non-zero disturbance bounds of  $\lvert w_j(t)\rvert \!\leq\! 0.02$, $j\!=\!1,2$, and $\lvert w_j(t)\rvert \!\leq\! 10^{-4}$, $j\!=\!1,2,3,4$, respectively, for all $t$. \textcolor{black}{All other RHOCP parameters remain the same for Examples~1 and~2 except in the former $\vert x_1(t) \vert \leq 20$, $j\!=\!1,2$, $x(0)\!= \![6.2\; 10]^\intercal$ and in the latter the horizon is shortened to $N\!=\!10$, $Q\! =\! 10^{-4}I$, $R \!=\! 10^{-4}I$, 
and the constraints are relaxed to  $\lvert x_j(t) \rvert \!\leq\! 20$, $j\!=\!1,2,3,4$ and $\lvert u_j(t) \rvert \leq 250$, $j\!=\!1,2$ $\forall t$. The initial seed in both examples was obtained using the initialization strategy described in Section~\ref{section:robustsection}.}
%
%
Figures~\ref{fig:trajplot_R} and \ref{fig:RcMPCplot1} show the closed loop trajectories for RcMPC applied to Example~1 with state $x_2$ reaching a minimum value of $-19.98$, thus robustly satisfying the constraints. 
Figure~\ref{fig:RcMPCplot1} shows that, in agreement with Theorem~\ref{thm:robustrecurfeas}, the seed tube $\mathbf{X}_1^{0,1}$ computed at the first iteration at time $i\!=\!1$ is contained within the optimal tube $\mathbf{X}_0^{\ast,3}$ computed at the final iteration for $i\!=\!0$.
Figure~\ref{fig:costplot_R1} demonstrates that, for any time step $i$, the optimal cost decreases at successive iterations and the perturbation sequence $\mathbf{c}_i^\ast$ converges to a value close to zero, as expected from Lemma~\ref{lem:JRcmpc} and Remark~\ref{rem:feastubesolut}. \textcolor{black}{Although the control perturbation converges to zero during the iteration at time step $i$, there is no guarantee that it will be smaller at first iteration at $i+1$. Thus, $\|\mathbf{c}_i^\ast\| \not\to 0$ as $i\!\to\!\infty$ in Figure~\ref{fig:costplot_R1}.}


\section{CONCLUSIONS}
This paper proposes a robust nonlinear model predictive control strategy for systems with convex state and control constraints, and with dynamics that can be expressed as a difference of convex functions. At each discrete time step, a sequence of convex problems is obtained by linearizing the system model around predicted trajectories. Robust tubes are used to bound model trajectories, thus accounting for linearization errors and external disturbances. In the absence of external disturbances we show that the successive approximation approach is non-conservative, since the predicted tubes and control sequences converge to a single trajectory which is a local optimum for the original nonlinear optimal control problem. We then reformulate the approach to consider additive disturbances, and derive a robust receding horizon strategy using a novel strategy for the selection of linearization points that provides a guarantee of recursive feasibility. We also establish a form of quadratic stability for the closed-loop system. As in the disturbance-free case, the robust approach solves a sequence of convex programming problems at each discrete time step without needing pre-computed error bounds or heuristics.
This iteration can be terminated early without affecting stability or feasibility guarantees.
\bibliographystyle{ieeetr}

\footnotesize{\bibliography{bibliography_edited}}

\section*{APPENDIX}
\normalsize\textit{Proof of Theorem \ref{thm:KKT}:} Theorem \ref{thm:convergence} shows that the cMPC iteration converges to a fixed point. We now prove that this point is a local optimum of the nMPC RHOCP. For readability we set $n_x=1$, but the proof for $n_x > 1$ is analogous. We also assume that the nMPC RHOCP \eqref{eqn:nMPC} satisfies a suitable constraint qualification (e.g.~linear independence constraint qualification, \cite{Reading72} Def.~12.4), so that a local minimizer satisfies the first order necessary conditions (\cite{Reading72} Thm.~12.1).

From Theorem \ref{thm:convergence} it is known that the cMPC iteration converges to $\{c_{ 0 \vert i},\ldots,c_{ N-1 \vert i}\}=\{0,\ldots,0\}$, $ S_{0 \vert i} = \{0\}$, \dots, $S_{N \vert i} = \{0\}$
and at this point we have
$(\mathbf{x}_i^\ast, \mathbf{u}_i^\ast) = (\mathbf{x}_i^{0,j},\mathbf{u}_i^{0,j})$ due to \eqref{eqn:pertrmodel}, \eqref{eqn:ueqn}. To prove that this point coincides with a local optimum for both cMPC and nMPC RHOCPs, we consider the equivalence to first-order of the cMPC RHOCP at this point and the nMPC RHOCP at $(\mathbf{x}_i ,\mathbf{u}_i) = (\mathbf{x}_i^{0,j},\mathbf{u}_i^{0,j}$). By inspection of \eqref{eqn:constr2}-\eqref{eqn:nmpccost} and \eqref{eqn:cMPCRHOCP2}-\eqref{eqn:cMPCRHOCP5} it is clear that if only constraints \eqref{eqn:constr2}, \eqref{eqn:constr3}, \eqref{eqn:constr1}, \eqref{eqn:constr5} and 
\eqref{eqn:cMPCRHOCP2}-\eqref{eqn:cMPCRHOCP1}, 
\eqref{eqn:cMPCRHOCP5} were present in the nMPC and cMPC RHOCPs, then their first order optimality conditions would be identical 
at the point $\mathbf{x}_i = \mathbf{x}^{0}_i$, $\mathbf{u}_i = \mathbf{u}_i^{0}$ and  $\{c_{ 0 \vert i},\ldots,c_{ N-1 \vert i}\}=\{0,\ldots,0\}$, $S_{0 \vert i} = \{0\}$, \dots, $S_{N \vert i} = \{0\}$. Therefore we prove that the same is true when \eqref{eqn:constr4} and \eqref{eqn:cMPCRHOCP4} are included, by showing that the following two simplified problems are equivalent to first-order.
\begin{subequations}
\begin{alignat}{2}
& \quad \text{\makebox[0pt][l]{$\textrm{nMPC simplified:}
\hspace{1mm}   \mathbf{c}_i^*  \defeq 
\argmin_{\mathbf{x}_i,\mathbf{c}_i} J^c_{\textrm{nMPC}}(\mathbf{x}_i,\mathbf{c}_i)$}}
& &\hspace{100mm}
\nonumber\\
&  \quad \text{\makebox[0pt][l]{subject to, for $k = 0, \ldots , N-1, x_{k+1 \vert i}- \phi(x_{k \vert i},c_{k \vert i}) =0 $}} 
\nonumber
\end{alignat}
\end{subequations}
%
where $f(x_{k \vert i},u_{k \vert i})= \phi(x_{k \vert i},c_{k \vert i})$ for $u_{k \vert i} = Kx_{k \vert i}+c_{k \vert i}$ and $J^c_{\textrm{nMPC}}(\cdot,\cdot)=J_{\textrm{nMPC}}(\cdot,\cdot)$, and
\begin{align}
 \refstepcounter{equation}
& \textrm{cMPC simplified:} \;\;\;\;  (\mathbf{c}_i^*, \mathbf{S}_i^*) = \argmin_{\mathbf{c}_i, \mathbf{S}_i} J^c_{\textrm{cMPC}}(\mathbf{c}_i, \mathbf{S}_i , \mathbf{x}_i^{0} ,\mathbf{c}_i^{0}) \notag \\
  &  \textrm{subject to, for $k=0,\ldots,N-1$ :} \;\;\;\;\;\;\;\;\;\; \;\;\;\;\;\;\;\;\;\; \;\;\;\;\;\;\;\;\;\;     \label{eqn:new} \\
   &  \bar{s}_{k+1 \vert i} \geq \phi(\bar{s}_{k \vert i}+x^{0}_{k \vert i},c_{k \vert i}+u^{0}_{k \vert i}-Kx^{0}_{k \vert i})-\phi(x^{0}_{k \vert i},c^{0}_{k \vert i}) \notag \\ 
      &  \bar{s}_{k+1 \vert i} \geq \phi(\munderbar{s}_{k \vert i}+x^{0}_{k \vert i}, c_{k \vert i}+u^{0}_{k \vert i}-Kx^{0}_{k \vert i})-\phi(x^{0}_{k \vert i},c^{0}_{k \vert i}) \notag \\ 
    & \munderbar{s}_{k+1 \vert i} \leq \Phi_{k \vert i} \bar{s}_{k \vert i}+ B_{k \vert i} c_{k \vert i} 
          \notag  \\
    & \munderbar{s}_{k+1 \vert i} \leq \Phi_{k \vert i} \munderbar{s}_{k \vert i} + B_{k \vert i} c_{k \vert i}
          \notag  
\end{align}
with $u^{0}_{k \vert i}=Kx^{0}_{k \vert i}+c^{0}_{k \vert i}$, 
$J^c_{\textrm{cMPC}}(\cdot, \cdot, \cdot ,\cdot)=J_{\textrm{cMPC}}(\cdot, \cdot , \cdot ,\cdot)$, and%
\begin{align*}
&\phi(s+x^0,c+  u^0-Kx^0) = f(s+x^0,Ks+c+u^0) \\
& \tfrac{\partial \phi(x,c)}{\partial x}  = \tfrac{\partial f(x,Kx+c)}{\partial x} =  \tfrac{\partial f(x,u)}{\partial x}+ \tfrac{\partial f(x,u)}{\partial u} \tfrac{\partial u}{\partial x} 
= A\!+\! BK =\Phi 
\\
&  \tfrac{\partial \phi(x,c)}{\partial c} =  \tfrac{\partial f(x,Kx+c)}{\partial c} =  \tfrac{\partial f(x,u)}{\partial u} \tfrac{\partial u}{\partial c} = B. \end{align*}
Defining the Lagrangian for the simplified nMPC RHOCP as $L_{\textrm{nMPC}}(\mathbf{x}_i ,\mathbf{c}_i,\mathbf{\lambda}_i) = 
   J^c_{\textrm{nMPC}}(\mathbf{x}_{i} ,\mathbf{c}_{i}) - \sum_{k=0}^N  \lambda_{k \vert i}[x_{k+1 \vert i}- \phi(x_{k \vert i},c_{k \vert i})]$,
the Karush-Kuhn-Tucker (KKT) first-order conditions \cite{Reading72} become:   
\begin{subequations}
\begin{align}
& 0 =\tfrac{\partial J^c_{\textrm{nMPC}}(\mathbf{x}_i,\mathbf{c}_i)}{\partial x_{0 \vert i}}\vert_{(\mathbf{x}_{i}^{*,n},\mathbf{c}_{i}^{*,n})}  + \lambda_{0 \vert i}^{*,n} \tfrac{\partial \phi(x,c)}{\partial x} \vert_{(x_{0 \vert i}^{*,n},c_{0 \vert i}^{*,n})} \label{eqn:nMPCKKT1} \\
&\textrm{and for }k=0,\ldots,N-1 \notag \\
&  0 =\tfrac{\partial J^c_{\textrm{nMPC}}(\mathbf{x}_i,\mathbf{c}_i)}{\partial x_{k \vert i}}\vert_{(\mathbf{x}_{i}^*,\mathbf{c}_{i}^{*,n})}  + \lambda_{k \vert i}^{*,n} \tfrac{\partial \phi(x,c)}{\partial x} \vert_{(x_{k \vert i}^{*,n},c_{k \vert i}^{*,n})}-\lambda^{*,n}_{k-1 \vert i} 
\label{eqn:nMPCKKT2}  \\
 & 0 = \tfrac{\partial J^c_{\textrm{nMPC}}(\mathbf{x}_i,\mathbf{u}_i)}{\partial u_{k \vert i}}\vert_{(\mathbf{x}_{i}^{*,n},\mathbf{c}_{i}^{*,n})} + \lambda_{k \vert i}^{*,n} \tfrac{\partial \phi(x,c)}{\partial c} \vert_{(x_{k \vert i}^{*,n},u_{k \vert i}^{*,n})}
\label{eqn:nMPCKKT3}  \\
& 0= x^{*,n}_{k+1 \vert i}- \phi(x^{*,n}_{k \vert i},c^{*,n}_{k \vert i}) \label{eqn:nMPCKKT4} 
\end{align}
\end{subequations}
where $(\cdot)^{*,n}$ denotes a local optimum point for nMPC simplified and $\mathbf{x}^{*,n}_i$ is obtained from $\mathbf{c}^{*,n}_i$ and $x_{k+1 \vert i}= \phi(x_{k \vert i},c_{k \vert i}) $.
For the simplified cMPC RHOCP we define the Lagrangian {\small $L_{\textrm{cMPC}}(\mathbf{c}_i,\mathbf{S}_i,\mathbf{x}^{0}_i,\mathbf{u}^{0}_i,\mathbf{\alpha}_i,\mathbf{\beta}_i,\mathbf{\gamma}_i,\mathbf{\mu}_i)=  J_{\textrm{cMPC}} \mathbf{c}_i,\mathbf{S}_i,\mathbf{x}^{0}_i,\mathbf{u}^{0}_i)
-\sum_{k=0}^{N-1}\alpha_{k \vert i}[\bar{s}_{k+1 \vert i}-\phi(\bar{z}_{k \vert i})+\phi(x^{0}_{k \vert i},c^{0}_{k \vert i})]
-\sum_{k=0}^{N-1}\beta_{k \vert i}[\bar{s}_{k+1 \vert i}  - \phi(\munderbar{z}_{k \vert i}) + \phi(x^{0}_{k \vert i},c^{0}_{k \vert i})]
-\sum_{k=0}^{N-1} \gamma_{k \vert i}(\Phi_{k \vert i} \bar{s}_{k \vert i} + B_{k \vert i} c_{k \vert i} -\munderbar{s}_{k+1 \vert i} )
-\sum_{k=0}^{N-1}\mu_{k \vert i}(\Phi_{k \vert i} \munderbar{s}_{k \vert i} + B_{k \vert i} c_{k \vert i}-\munderbar{s}_{k+1 \vert i} )$} with KKT conditions:%
\begin{subequations}%
\begin{align} 
  &0 =  \tfrac{\partial J_{\textrm{cMPC}}  }{\partial \bar{s}_{0 \vert i}}\vert_{\mathbf{z}_i^{*,c}} + \alpha_{0 \vert i}^{*,c} \tfrac{\partial \phi(x,c)}{\partial x} \vert_{\bar{z}^{*,c}_{0 \vert i}} -\gamma_{0 \vert i}^{*,c} \tfrac{\partial \phi(x^{0}_{0 \vert i},c^{0}_{0 \vert i})}{\partial x^{0}_{0 \vert i}} \label{eqn:cMPCKKT1} \\
  &0=  \tfrac{\partial J_{\textrm{cMPC}}  }{\partial \munderbar{s}_{0 \vert i}}\vert_{\mathbf{z}_i^{*,c}} +\beta_{0 \vert i}^{*,c} \tfrac{\partial \phi(x,c)}{\partial x} \vert_{\munderbar{z}^{*,c}_{0 \vert i}}- \mu_{k \vert i}^{*,c} \tfrac{\partial \phi(x^{0}_{k \vert i},c^{0}_{k \vert i})}{\partial x^{0}_{k \vert i}} \label{eqn:cMPCKKT2} \\
&\textrm{and for } k =1,\ldots,N-1\notag \\
&0=  \tfrac{\partial J_{\textrm{cMPC}}  }{\partial \bar{s}_{k \vert i}}\vert_{\mathbf{z}^{*,c}_i} -(\alpha_{k-1 \vert i}^{*,c}+ \beta^{*,c}_{k-1 \vert i})
+ \alpha_{k \vert i}^{*,c} \tfrac{\partial \phi(x,u)}{\partial x} \vert_{\bar{z}^{*,c}_{k \vert i} } 
\nonumber \\
&\qquad -\gamma_{k \vert i}^{*,c} \tfrac{\partial \phi(x^{0}_{k \vert i},c^{0}_{k \vert i})}{\partial x^{0}_{k \vert i}} \label{eqn:cMPCKKT3}\\
&0=  \tfrac{\partial J_{\textrm{cMPC}}  }{\partial \munderbar{s}_{k \vert i}}\vert_{\mathbf{z}_i^{*,c}} +(\mu_{k-1 \vert i}^{*,c}+\gamma_{k-1 \vert i}^{*,c}) +\beta_{k \vert i}^{*,c} \tfrac{\partial \phi(x,c)}{\partial x} \vert_{\munderbar{z}^{*,c}_{k \vert i} } \nonumber \\
&\qquad - \mu_{k \vert i}^{*,c} \tfrac{\partial \phi(x^{0}_{k \vert i},c^{0}_{k \vert i})}{\partial x^{0}_{k \vert i}} \label{eqn:cMPCKKT4} \\
&0=  \tfrac{\partial J_{\textrm{cMPC}} }{\partial c_{k \vert i}}\vert_{\mathbf{z}^{*,c}_i} + \alpha_{k \vert i}^{*,c} \tfrac{\partial \phi(x,c)}{\partial c} \vert_{\bar{z}^{*,c}_{k \vert i} }+\beta_{k \vert i}^{*,c} \tfrac{\partial \phi(x,c)}{\partial c} \vert_{\munderbar{z}^{*,c}_{k \vert i} }
\nonumber \\
&\qquad - (\gamma_{k \vert i}^{*,c}+\mu_{k \vert i}^{*,c}) \tfrac{\partial \phi(x^{0}_{k \vert i},c^{0}_{k \vert i})}{\partial c^{0}_{k \vert i}} \label{eqn:cMPCKKT5} \\
&0=\alpha_{k \vert i}(\bar{s}^{*,c}_{k+1 \vert i} -\phi(\bar{z}^{*,c}_{k \vert i} )+\phi(x^{0}_{k \vert i},u^{0}_{k \vert i})) \label{eqn:cMPCKKT7} \\
&0= \beta_{k \vert i}(\bar{s}^{*,c}_{k+1 \vert i} -\phi(\munderbar{z}^{*,c}_{k \vert i} )+\phi(x^{0}_{k \vert i},u^{0}_{k \vert i}))  \label{eqn:cMPCKKT8} \\
&0= \gamma_{k \vert i}(\Phi_{k \vert i} \bar{s}^{*,c}_{k \vert i} + B_{k \vert i} c^{*,c}_{k \vert i} -\munderbar{s}^{*,c}_{k+1 \vert i} ) \label{eqn:cMPCKKT9} \\
&0=\mu_{k \vert i}(\Phi_{k \vert i} \munderbar{s}^{*,c}_{k \vert i} + B_{k \vert i} c^{*,c}_{k \vert i}-\munderbar{s}^{*,c}_{k+1 \vert i} ) \label{eqn:cMPCKKT10} \\
& \alpha^{*,c}_{k \vert i}\geq0,\;\; \beta^{*,c}_{k \vert i} \geq 0, \;\; \gamma^{*,c}_{k \vert i} \geq 0, \;\; \mu^{*,c}_{k \vert i} \geq 0 \label{eqn:cMPCKKT11}
\end{align}
\end{subequations}
where $\mathbf{z}_{i}^{*,c} = (\mathbf{c}_i^{*,c},\mathbf{S}_i^{*,c},\mathbf{x}^{0}_i,\mathbf{u}^{0}_i)$ and 
\begin{align*}
\nabla_{c_{k \vert i}}\phi (\munderbar{s}_{k \vert i}+x^{0}_{k \vert i},c_{k \vert i}+u^{0}_{k \vert i}-Kx^{0}_{k \vert i}) 
&= \nabla_{c}\phi(x,c)\vert_{\munderbar{z}_{k \vert i}}  \\
\nabla_{\munderbar{s}_{k \vert i}}\phi (\munderbar{s}_{k \vert i}+x^{0}_{k \vert i},c_{k \vert i}+u^{0}_{k \vert i}-Kx^{0}_{k \vert i}) 
&= 
\nabla_{x}\phi(x,u) \vert_{\munderbar{z}_{k \vert i}}
\\
(\bar{s}_{k \vert i}+x^{0}_{k \vert i},c_{k \vert i}+u^{0}_{k \vert i}-Kx^{0}_{k \vert i})
&=
\bar{z}_{k \vert i} 
\\
(\munderbar{s}_{k \vert i}+x^{0}_{k \vert i},c_{k \vert i}+u^{0}_{k \vert i}- Kx^{0}_{k \vert i} ) 
&= 
\munderbar{z}_{k \vert i} 
\\
(\bar{s}^{*,c}_{k \vert i}+x^{0}_{k \vert i},c^{*,c}_{k \vert i}+u^{0}_{k \vert i} - Kx^{0}_{k \vert i}) 
&=\bar{z}^{*,c}_{k \vert i} 
\\
(\munderbar{s}^{*,c}_{k \vert i}+x^{0}_{k \vert i},c^{*,c}_{k \vert i}+u^{0}_{k \vert i} - Kx^{0}_{k \vert i}) 
&=
\munderbar{z}^{*,c}_{k \vert i} 
\end{align*}
where $(\cdot)^{*,c}$ denotes solution of cMPC simplified and  $\mathbf{x}^{*,c}_i$ is obtained from $\mathbf{c}^{*,c}_i$ and $x_{k+1 \vert i}= \phi(x_{k \vert i},c_{k \vert i}) $.
For each $k$, at least one of the constraints in \eqref{eqn:new} must be active due to the definition of $J_{\textrm{cMPC}}^c$, but all four constraints in \eqref{eqn:new} cannot be simultaneously active (unless $S_{k \vert i}=\{0\}$, in which case $\bar{s}_{k \vert i}=\munderbar{s}_{k \vert i}$),
so at least one multiplier is non-zero and at least one is zero. 
Therefore either $ \tfrac{\partial J_{\textrm{cMPC}}  }{\partial \bar{s}_{k \vert i}}\vert_{\mathbf{z}^{*,c}_i}$ or $ \tfrac{\partial J_{\textrm{cMPC}}  }{\partial \munderbar{s}_{k \vert i}}\vert_{\mathbf{z}^{*,c}_i}$ is zero in \eqref{eqn:cMPCKKT3}, \eqref{eqn:cMPCKKT4} for each $k$;
of these, the nonzero gradient is equal to $ \tfrac{\partial J_{\textrm{nMPC}}  }{\partial x_{k \vert i}}\vert_{(\mathbf{x}^r_i,\mathbf{c}^{r}_i)}$ due to the definition of $J_{\textrm{cMPC}}$,
where $c^r_{k \vert i} = c^{*,c}_{k \vert i}+c^{0}_{k \vert i}$ and 
$x^r_{k \vert i} = \bar{s}_{k \vert i}^{*,c}+x^{0}_{k \vert i}$ or $\munderbar{s}_{k \vert i}^{*,c}+x^{0}_{k \vert i}$ depending on which constraint in \eqref{eqn:new} is active. 
From \eqref{eqn:cMPCKKT1}-\eqref{eqn:cMPCKKT4} we have
\begin{align}          0&=  \tfrac{\partial J_{\textrm{nMPC}}  }{\partial x_{0 \vert i}}\vert_{(\mathbf{x}^r_i,\mathbf{c}^{r}_i)} + \alpha_{0 \vert i}^{*,c} \tfrac{\partial \phi(x,c)}{\partial x} \vert_{\bar{z}^{*,c}_{0 \vert i} }  +\beta_{0 \vert i}^{*,c} \tfrac{\partial \phi(x,c)}{\partial x} \vert_{\munderbar{z}^{*,c}_{0 \vert i} } \notag \\
     & \;\;\;\; - (\gamma_{0 \vert i}^{*,c} +\mu_{0 \vert i}^{*,c}) \tfrac{\partial \phi(x^{0}_{0 \vert i},c^{0}_{0 \vert i})}{\partial x^{0}_{0 \vert i}} \label{eqn:cMPCx0grad} \\
          0&=  \tfrac{\partial J_{\textrm{nMPC}}  }{\partial x_{k \vert i}}\vert_{(\mathbf{x}^r_i,\mathbf{c}^{r}_i)} + \alpha_{k \vert i}^{*,c} \tfrac{\partial \phi(x,c)}{\partial x} \vert_{\bar{z}^{*,c}_{k \vert i} }  +\beta_{k \vert i}^{*,c} \tfrac{\partial \phi(x,c)}{\partial x} \vert_{\munderbar{z}^{*,c}_{k \vert i} } \notag \\
     & \;\;\;\; - (\gamma_{k \vert i}^{*,c} +\mu_{k \vert i}^{*,c}) \tfrac{\partial \phi(x^{0}_{k \vert i},c^{0}_{k \vert i})}{\partial x^{0}_{k \vert i}} \notag \\& \;\;\;\; -(\alpha_{k-1 \vert i}^{*,c}+ \beta^{*,c}_{k-1 \vert i}-\mu_{k-1 \vert i}^{*,c}-\gamma_{k-1 \vert i}^{*,c}) . \label{eqn:cMPCxkgrad}
\end{align}
By Theorem \ref{thm:convergence}, the cMPC iteration converges to a fixed point with $\mathbf{c}^\ast_i=0$, $\mathbf{S}^\ast_i = \{\{0\}, \ldots, \{0\}\}$, which we denote $(0,\{0\}, \mathbf{x}^{0}_{i},\mathbf{u}^{0}_{i})$. As this point is optimal for the cMPC RHOCP, the simplified cMPC KKT conditions must be satisfied. Note that 
$(\mathbf{x}^{r}_{i}, \mathbf{c}^{r}_{ i}) = (\mathbf{x}^{0}_{ i}, \mathbf{c}^{0}_{i})$ if $\mathbf{S}^\ast_{i } = \{\{0\},\ldots,\{0\}\}$ and $\mathbf{c}^\ast_{i}=0$; hence comparing \eqref{eqn:nMPCKKT1} with \eqref{eqn:cMPCx0grad}, \eqref{eqn:nMPCKKT2} with \eqref{eqn:cMPCxkgrad}), and \eqref{eqn:nMPCKKT3} with \eqref{eqn:cMPCKKT5} for $(\mathbf{x}^{*,n}_{i},\mathbf{c}^{*,n}_{i})=(\mathbf{x}^{0}_{i},\mathbf{c}^{0}_{i})$, it can be seen that they are equivalent  if
     $\lambda^{*,n}_{k \vert i} =  \alpha^{*,c}_{k \vert i}+\beta^{*,c}_{k \vert i}-\gamma^{*,c}_{k \vert i}-\mu^{*,c}_{k \vert i}$ for all $k$.
Furthermore Theorem \ref{thm:recrsfeas} ensures \eqref{eqn:nMPCKKT4} is satisfied by $(\mathbf{x}^{*,n}_{i},\mathbf{c}^{*,n}_{i}) = (\mathbf{x}^{0}_{i},\mathbf{c}^{0}_{i})$. 
Hence the KKT conditions for cMPC and nMPC RHOCPs are equivalent at the cMPC optimum.
\qed 

\begin{IEEEbiography}[{\includegraphics[width=1in,height=1.25in,clip,keepaspectratio]{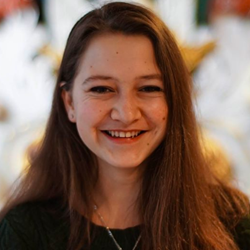}}]{Yana Lishkova} completed a Bachelor and MEng program at the University of Cambridge focusing on Aerospace Engineering and Instrumentation \& Control. She is currently a D.Phil. student in the Control Group at Oxford University. Her work investigates novel methods for simulation, optimization and control of multirate systems, which exhibit dynamics on different time scales. 
Applications include safety-critical problems in the aerospace domain. Her research is supported by an EPSRC Research Studentship. In 2022 she was awarded the Amelia Earhart Fellowship. 
\end{IEEEbiography}

\begin{IEEEbiography}[{\includegraphics[width=1in,height=1.25in,clip,keepaspectratio]{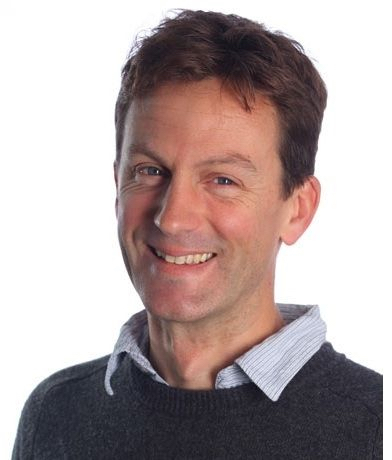}}]{Mark Cannon} received the M.Eng.~and D.Phil.~degrees from the University of Oxford in 1993 and 1998, respectively, and the M.S.~degree from the Massachusetts Institute of Technology in 1995. He is currently an Associate Professor of Engineering Science at Oxford University and Fellow of St. John’s College, Oxford. His research interests include 
robust and stochastic receding horizon control, learning and optimization-based control, data-driven control and optimization, adaptive control, and applications in aerospace, automotive and biomedical engineering.
\end{IEEEbiography}

\end{document}